\documentclass[oneside]{amsart}
\usepackage{amsmath,amssymb,hyperref,mathrsfs,graphicx,bbm}
\usepackage[utf8x]{inputenc}
\usepackage{amsmath}
\usepackage{amsfonts}
\usepackage{latexsym}
\usepackage{amsthm}
\usepackage{stmaryrd,amscd}
\usepackage{xargs}
\usepackage{tikz}
\usepackage{tikz-cd}
\usepackage{graphicx}
\usepackage{subcaption}
\usepackage{color}
\usepackage{enumitem}

\usepackage[colorinlistoftodos,prependcaption]{todonotes}
\presetkeys{todonotes}{inline,backgroundcolor=white}{}
\tikzset{/tikz/notestyleraw/.append style={text=black}}

\newtheorem{thm}{Theorem}[section]

\newtheorem{lem}[thm]{Lemma}
\newtheorem{defn}[thm]{Definition}
\newtheorem{prop}[thm]{Proposition}
\newtheorem{cor}[thm]{Corollary}

\newtheorem{que}{Question}
\newtheorem{rmk}[thm]{Remark}

\numberwithin{equation}{section}

\newcommand{\be}{\begin{eqnarray}}
\newcommand{\ee}{\end{eqnarray}}
\newcommand{\beal}{\begin{aligned}}
\newcommand{\enal}{\end{aligned}}
\newcommand{\ben}{\begin{eqnarray*}}
\newcommand{\een}{\end{eqnarray*}}
\newcommand{\beq}{\begin{equation}}
\newcommand{\eeq}{\end{equation}}

\newcommand{\ess}{\operatorname{ess}}

\newcommand{\eps}{\epsilon}

\newcommand{\lb}{\lambda}
\newcommand{\dt}{\delta}

\newcommand{\om}{\omega}
\newcommand{\Om}{\Omega}

\newcommand{\T}{\mathbb{T}}
\newcommand{\N}{\mathbb{N}}
\newcommand{\R}{\mathbb{R}}

\newcommand{\Z}{\mathbb{Z}}

\newcommand{\cD}{\mathcal{D}}

\newcommand{\cA}{\mathcal{A}}
\newcommand{\cL}{\mathcal{L}}

\newcommand{\cT}{\mathcal{T}}

\newcommand{\wt}{\widetilde }

\title{The breakdown of invariant graphs for dissipative standard maps}

\thanks{{\it Email address:} jellychung1987@gmail.com$^\dagger$, zhourui@amss.ac.cn$^\ddagger$}
\begin{document}
\maketitle

\centerline{Jianlu Zhang$^\dagger$, Rui Zhou$^\ddagger$}
\medskip
{\footnotesize
\centerline{State Key Laboratory of Mathematical Sciences,}
 \centerline{Academy of Mathematics and systems science,}
 \centerline{Chinese Academy of Sciences, Beijing 100190, China.}
}
\medskip

\begin{abstract}

We are concerned with the dissipative standard maps $\Phi: T^*\T\rightarrow T^*\T$ of the form
\[
\left\{
\begin{aligned}
&X=x+\frac{1-e^{-\lb}}{\lb} y+\om,\\
&Y=e^{-\lb}y+\eps e^{-\lb}V'(X),
\end{aligned}
\right.
\]
with $\lb, \eps>0$, $V(\cdot)\in C^2(\R\slash\Z,\R)$ and $\om\in\R$, focusing mainly on the mechanism by which KAM tori (i.e. KAM-type invariant graphs) get twisted and become non-graphic as $\eps$ rises. We first look into the variational properties of invariant graphs, then establish computational upper and lower bounds for the breakdown threshold, together with a numerical algorithm that evaluates them. We also discuss the mechanism underlying the breakdown and the drift of rotation numbers for invariant graphs as $\eps$ tends to the threshold.
\end{abstract}

\section{\label{s1}Introduction}

The complex structure of global attractors for dissipative twist maps was studied in \cite{Cal,leCal2}. A basic prototype is the DSM, for which one can track how the simple attractors---unique invariant graphs---get twisted as $\eps$ increases. This is a direct analogue of KAM-torus destruction for area-preserving standard maps \cite{Mat1,Tre,Mo}, but the dissipative dynamics differs enough to require separate treatment.

In recent works \cite{Ar,CCD,CCD1,MS}, the concept {\it conformally symplectic twist map} was introduced and refers to a large class of dissipative twist maps
(including dissipative standard maps) possessing variational advantages. Precisely, one can always find a certain {\it Aubry--Mather set} lying within the global attractor,
which in some sense imposes additional variational constraints on the structure of attractors. Such a specific set is always graphic \cite{Cal, MS}, on which each trajectory contains no conjugate points. In particular, if an invariant graph is of KAM type, say the map restricted to it is conjugated to a rotation, then it has to be {\it normally hyperbolic} and therefore is the only invariant set in its suitably small neighbourhood, see \cite{CCD} for the proof. Such a feature enables the persistence of KAM graphs under $C^1-$perturbations of the map,
whereas the drift of rotation numbers under perturbations is totally unpredictable.

In this paper, we investigate the breakdown threshold of invariant graphs for dissipative standard maps and develop a numerical scheme to estimate the critical values of $\eps$ from both the upper and lower sides. We also look into the mechanism of the breakdown of invariant graphs and analyse the drift of rotation numbers for the Aubry--Mather set, which always persists as the broken part of the invariant graphs as $\eps$ increases. Unlike \cite{CC}, where Newton's method is used (so that a drifting parameter is introduced to maintain the Diophantine rotation number), our approach relies mainly on variational methods. We also examine the drift of rotation numbers and identify possible mechanisms therein, by analysing the structural catastrophe of the global attractor, which locks the rotation numbers as $\eps$ reaches the threshold.

We finally mention \cite{SW}, which proves breakdown of invariant graphs for integrable CSTMs under low-regularity ($C^{1-\eps}$) perturbations. Here we consider large, higher-regularity perturbations and discuss the mechanism.

\medskip

\noindent{\bf Notation and regularity.} Throughout, $V$ is the periodic potential and $V'$ the forcing term; both are periodic and $\int_\T V'=0$. For the DSM we assume $V\in C^2(\T,\R)$ (equivalently $V'\in C^1$) unless stated otherwise. In the generic mechanism discussion of Section~\ref{s4} we also assume $V'$ has a unique maximum and minimum; results not using this are indicated. The smooth-linearization argument in Subsection~\ref{subsec:saddle-node} needs stronger smoothness and non-resonance, specified there.\medskip

\noindent{\bf Terminology.} A {\it circle} is a topologically non-contractible closed curve in $T^*\T$, and an {\it invariant circle} is a circle that is $\Phi$-invariant (equivalently, an invariant subset of $T^*\T$ homeomorphic to $\T$). An invariant circle is {\it graphic} if it is the graph of a continuous function $y=\psi(x)$, $\psi:\T\rightarrow\R$, and {\it non-graphic} otherwise; a graphic invariant circle is also called an {\it invariant graph}. The {\it breakdown} of an invariant graph means the loss of the graphic property of a given invariant circle as a parameter varies.\medskip

\noindent{\bf Parameter conventions.} Throughout, $\lambda>0$ is fixed. For the DSM we fix $\lambda,\omega,V$ and vary $\eps\ge0$; $\eps_*(\lambda,\omega)$ is defined for that family. To study its dependence on $\omega$ we compare families at the same $\lambda$. Numerical experiments use $\lambda=2$ unless stated otherwise.\medskip

\noindent{\bf Structure of the paper.} This article is organised as follows:
\begin{itemize}
\item In Section \ref{s1.1} we formulate the variational framework for conformally symplectic twist maps (CSTM) and reveal the weak KAM properties associated.
    \item In Section~\ref{s2}, we prove that any invariant graph has to be globally minimal, in the sense that it coincides with the Aubry set and the gradient graph of a unique discrete weak KAM solution, see Theorem \ref{thm:lip} and Theorem \ref{thm:aubry}.
    \item In Section~\ref{s3}, we consider the breakdown of invariant graphs for a specific kind of CSTMs, i.e. the dissipative standard maps. We obtain an upper bound and a lower bound estimate for the threshold of the breakdown of invariant graphs, see Theorem~\ref{thm:K_H} and Theorem \ref{thm:e_u_i}. Both bounds are computable and numerical evidence is also supplied.
    \item In Section~\ref{s4}, we discuss the mechanism of the breakdown of invariant graphs by summarizing all possible structural catastrophes. We also exhibit numerical results and figures to show some new phenomena, such as the locking of the drift and gaps in the rotation numbers of the Aubry--Mather set.
\item In Appendix~\ref{app:error}, we describe the numerical algorithm used to evaluate the finite-iteration thresholds, and quantify its truncation and discretisation errors.
\end{itemize}

\medskip

\noindent{\bf Acknowledgements.} This work is supported by the National Key R\&D Program of China (No. 2022YFA1007500) and the National NSF of China (No. 12231010, 12571207). The authors are grateful to Lin Wang for helpful discussions and suggestions. We wish to acknowledge the support of the author community in using REV\TeX{}, offering suggestions and encouragement, testing new versions. All data and figures in this paper were produced with Python, using NumPy, SciPy and Matplotlib, with Pillow for composing the figures, pandas and openpyxl for handling the data tables, and numba, mpmath and tqdm in some auxiliary scripts.

\medskip

\noindent{\bf Statements and Declarations.} The authors declare no competing interests and no data associated.

\medskip

\noindent{\bf Generative AI disclosure.} LLM tools were used to assist with parts of the programming and with grammar and spelling checks.

\section{Conformally symplectic twist map and the weak KAM properties}\label{s1.1}

This section concerns general exact conformally symplectic twist maps. We assume that the map is a $C^1$ diffeomorphism and that its generating function is $C^2$, with the twist, periodicity, and coercivity properties specified below. Throughout, $\T:=\R/\Z$ denotes the flat circle and $T^*\T:=\T\times\R$ the cylinder. A function on $\T$ is identified with its $1$-periodic lift to $\R$; thus, when $x\in\R$ appears as the argument of such a function, its value is understood modulo $\Z$.

A $C^1$-diffeomorphism $\Phi: T^*\T \rightarrow T^*\T$ formally expressed by
\beq\label{eq:cstm}
\Phi:
\begin{cases}
X = f(x,y),\\
Y = g(x,y),
\end{cases}
\eeq
is called a \textit{conformally symplectic twist map (CSTM)} if it satisfies:
\begin{itemize}
    \item[(a)] $dX\wedge dY=e^{-\lb}dx\wedge dy$ for some $\lb>0$;
    \item[(b)] $\lim_{|y|\rightarrow+\infty}|X-x|=+\infty, \forall x\in\T$;
    \item[(c)] $\partial_y f>0$ for all $(x,y)\in T^*\T$;
\end{itemize}
Furthermore, if we can find a $C^2$ smooth function $h: \R \times \R \to \R$ satisfying
\begin{itemize}
    \item[(d)] $h(x+1, X+1) = h(x,X)$ for all $(x,X) \in \R^2$;
    \item[(e)] $\lim_{|X-x| \to +\infty} h(x,X) = +\infty, \forall x\in\T$;
    \item[(f)] $-\partial_{12} h(x,X) > 0$ everywhere on $\R^2$;
\end{itemize}
such that
\beq\label{eq:gene-map}
y = -\partial_1 h(x,X), \qquad Y = e^{-\lb} \partial_2 h(x,X),
\eeq
then the associated $\Phi$ will be called {\it exact}. For any exact CSTM $\Phi$, the {\it flux} number

\beq
\operatorname{Flux}(\Phi) := \oint_z^{Tz} e^\lb Y \, dX - y \, dx,
\eeq
is always zero whatever $z=(x,y)\in\R^2$ is chosen ($Tz:=(x+1,y)\in\R^2$) due to the {\it Stokes' Theorem}. Such an intrinsic conservation supplies fine variational properties of $\Phi$. To stress this, we need the following definition:

\begin{defn}{\label{defn:crit_conf}}
A bidirectional sequence $\{\bar x_i \in \mathbb{R}\}_{i \in \mathbb{Z}}$ is called a \textit{\bf critical configuration}, if
\be\label{eq:disc-e-l}
\partial_1 h(\bar x_i, \bar x_{i+1}) + e^{-\lambda} \partial_2 h(\bar x_{i-1}, \bar x_i) = 0
\ee
for all $i \in \mathbb{Z}$.
Any critical configuration determines a trajectory $\left\{\left(x_i, -\partial_1 h(x_i, x_{i+1})\right)\right\}_{i \in \mathbb{Z}}$ of \eqref{eq:cstm}, and vice versa\footnote{Throughout this paper, we put `$-$' on $x_i$ to represent its lift to the universal space $\R$. In most of the occasions, we can treat $[0,1)\subset\R$ equally with $\T$.}. We call a critical configuration $\{x_i \in \mathbb{R}\}_{i \in \mathbb{Z}}$ {\bf fixed-endpoint minimizing}, if for all integers $n\leq m$, the {\bf action function}
\be\label{eq:act}
\Lambda_n^m(x_n, x_m) := \inf_{\substack{\{s_i\}_{i=n}^m\subset\R \\ s_n = x_n,\ s_m = x_m}} \sum_{i=n}^{m-1} e^{\lambda i} h(s_i, s_{i+1}).
\ee
achieves the minimum via $s_i=x_i$ with $n\leq i\leq m$.
\end{defn}

\subsection{Discrete weak KAM solution}

The weak KAM theory was originally used to study the Aubry--Mather sets of positive definite Hamiltonian systems, and was recently generalized to conformally symplectic Hamiltonians in \cite{MS, WWY}. We briefly make a discrete adaptation of this theory to the current settings, see also \cite{DFIZ,SW} for some relevant works.

For any $n\in\N$, the {\bf backward Lax--Oleinik operator} $\cL^-_n :C(\T,\R)\rightarrow C(\T,\R)$ is defined by
\be\label{eq:l-o}
\cL_n^-\psi(x)=\inf_{q\in\R}\big\{\, e^{-\lb n}\psi(q)+\Lambda_{-n}^{0}(q,x)\,\big\}.
\ee
As a convention we always take $\cL_0^-=id$.
\begin{prop}\label{prop:weak-kam}
There exists a unique fixed point $u\in {\rm Lip}(\T,\R)$ of $\cL_n^-$, i.e. $\cL_n^-u(x)=u(x)$ for all $n\in\N$. Such a $u(x)$ is called a {\bf (discrete) weak KAM solution} of \eqref{eq:cstm}, in the sense that
\begin{itemize}
\item[(i)] $
e^{\lb m}u(x_1)-e^{\lb n}u(x_0)\leq \Lambda_n^m(\bar x_0,\bar x_1)$ for any integers $n\leq m$;
\item[(ii)] for any $x\in\T$, there exists a backward {\bf calibrated} configuration $\{ x_i\}_{i=-\infty}^0$ with $x_0=x$, such that for any $n\leq m\leq 0$, there holds
\be\label{eq:cal}
e^{\lb m}u(x_m)-e^{\lb n}u(x_n)=\sum_{i=n}^{m-1} e^{\lb i} h(x_i,x_{i+1}).
\ee
\end{itemize}
\end{prop}
\begin{proof}
Fix $\psi\in C(\T,\R)$. By the coercivity (e) of $h$, the minimizing configurations of $\cL_n^-\psi$ satisfy, uniformly in $n$,
\beq\label{eq:compact}
\max_{-n+1\leq i\leq 0}|\bar x_i-\bar x_{i-1}|\leq K,
\eeq
with $K$ depending only on $\lb$ and $h$. This yields the uniform boundedness and Lipschitz continuity of $\{\cL_n^-\psi\}_{n\ge1}$, hence the existence of limit points; since $\cL_1^-$ is a contraction of ratio $e^{-\lb}$ in the sup norm, the limit $u$ is unique. The domination inequality (i) follows from the definition of $\cL_n^-$, and the calibration identity (ii) from the compactness \eqref{eq:compact}; the details are as in \cite{Fa,MS}.
\end{proof}

\begin{lem}
Any weak KAM solution $u$ is {\bf semiconcave} with quadratic modulus, i.e.
\[
e^{\lb}\bigl(u(x+\eta)+u(x-\eta)-2u(x)\bigr)\leq C|\eta|^2
\]
for any $x\in\T$ and $|\eta|\leq 1/4$.
\end{lem}

\proof
Let $\{\bar x_i\}_{i\le0}$ be the backward calibrating configuration with $\bar x_0=x$, so that $e^\lb u(x)=u(\bar x_{-1})+h(\bar x_{-1},x)$ and $e^\lb u(x\pm\eta)\le u(\bar x_{-1})+h(\bar x_{-1},x\pm\eta)$. Then
\ben
e^\lb\bigl(u(x+\eta)+u(x-\eta)-2u(x)\bigr)
&\le& h(\bar x_{-1},x+\eta)+h(\bar x_{-1},x-\eta)-2h(\bar x_{-1},x)\\
&\le&\max_{y\in[x-\eta,x+\eta]}|\partial_{22}h(\bar x_{-1},y)|\,|\eta|^2,
\een
the maximum being finite because $|y-\bar x_{-1}|\le c_1(\lb)$ by \eqref{eq:compact}.
\qed

\begin{prop}
A configuration $\{\bar x_i\}_{i\in\Z}$ is called {\bf globally minimal} (also called {\it calibrated}, see \cite{Fa}), if for any $n<m$, we have
\beq
e^{\lb m}u(x_m)-e^{\lb n}u(x_n)=\sum_{i=n}^{m-1} e^{\lb i} h(\bar x_i,\bar x_{i+1}).
\eeq
The following conclusions can be drawn:
\begin{itemize}
\item[(1)] There exists at least one globally minimal configuration;
\item[(2)] Each globally minimal configuration uniquely determines a trajectory of $\Phi$, i.e. for any $i\in\Z$, there exists a unique $y_i\in\R$ such that $\{(x_i,y_i)\}_{i\in\Z}$ is a trajectory of $\Phi$;
\item[(3)] $u$ is differentiable on $\{x_i\}_{i\in\Z}$;
\end{itemize}
\end{prop}
\proof
For $x\in\T$ let $\{\bar x_i\}_{i\in\Z_-}$ be a backward configuration ending with $x$ and satisfying \eqref{eq:cal}; it is minimal in the sense of \eqref{eq:act} for $n<m\le0$. Then $y_i:=e^{-\lb}\partial_2h(\bar x_{i-1},\bar x_i)$ defines a backward trajectory, uniformly bounded by \eqref{eq:compact}, and its $\alpha$-limit set $\Om_\alpha$ is invariant; by \eqref{eq:cal} every orbit in $\Om_\alpha$ determines a globally minimal configuration, which proves (1). For (2), a globally minimal configuration is fixed-endpoint minimizing, so $e^\lb u(x_{i+1})-e^{-\lb}u(x_{i-1})=e^{-\lb}h(\bar x_{i-1},\bar x_i)+h(\bar x_i,\bar x_{i+1})$ with $\bar x_i$ minimizing the right-hand side, and $y_i:=-\partial_1h(x_i,x_{i+1})$ is uniquely determined. For (3), applying the domination inequality of Proposition~\ref{prop:weak-kam}(i) to the segments before and after $x_{i+1}$ gives
\[
e^{-\lb}\bigl(h(x_i,x_{i+1})-h(x_i,y)\bigr)\ \le\ u(x_{i+1})-u(y)\ \le\ h(y,x_{i+2})-h(x_{i+1},x_{i+2}),
\]
so $du(x_{i+1})$ exists and equals $e^{-\lb}\partial_2h(x_i,x_{i+1})$.
\qed
\begin{defn}\label{def:aubry}
Denote by $\cA$ the union of the $\alpha$-limit sets of all backward calibrated configurations associated with the weak KAM solution $u$. By item (3) of the preceding proposition, $u$ is differentiable on $\cA$. Then we define the {\bf Aubry set} $\wt\cA:=\{(x,du(x))\in T^*\T\mid x\in\cA\}$, which is invariant under $\Phi$, and the canonical projection $\pi:\wt\cA\rightarrow\cA$ is a Lipschitz function.
\end{defn}

\begin{lem}\label{lem:aubry-fixed-endpoint}
Every orbit $\{(x_i,y_i)\}_{i\in\Z}$ contained in the Aubry set $\wt\cA$ is fixed-endpoint minimizing. More precisely, for any integers $n<m$ and any coherent choice of lifts $\{\bar x_i\}_{i=n}^m$ of the orbit, one has
\beq\label{eq:aubry-fixed-endpoint}
\sum_{i=n}^{m-1}e^{\lb i}h(\bar x_i,\bar x_{i+1})=\Lambda_n^m(\bar x_n,\bar x_m).
\eeq
\end{lem}

\begin{proof}
Let $z_0=(x_0,du(x_0))\in\wt\cA$. By Definition~\ref{def:aubry} it is an $\alpha$-limit point of a backward calibrated orbit $\{z_i^-\}_{i\le0}$; choose $k_j\to+\infty$ with $z^-_{-k_j}\to z_0$. Then $z^-_{i-k_j}\to z_i:=\Phi^i(z_0)$ for every fixed $i$, with lifts chosen consistently. For large $j$ the identity \eqref{eq:cal} holds between $n-k_j$ and $m-k_j$; multiplying it by $e^{\lb k_j}$ and letting $j\to\infty$ gives $e^{\lb m}u(x_m)-e^{\lb n}u(x_n)=\sum_{i=n}^{m-1}e^{\lb i}h(\bar x_i,\bar x_{i+1})$, while the domination inequality of Proposition~\ref{prop:weak-kam}(i) gives the reverse inequality for every competing configuration. Hence \eqref{eq:aubry-fixed-endpoint}.
\end{proof}

\begin{rmk}\label{rem:fixed-endpoint-converse}
The converse is false. For the DSM with $\eps=0$, the Aubry set is $\wt\cA=\T\times\{0\}$, but every orbit is fixed-endpoint minimizing. Indeed, for $s_i=\bar x_i+\delta_i$ with $\delta_n=\delta_m=0$, the
weighted action difference is
\[
\sum_{i=n}^{m-1}e^{\lb i}\bigl[h(s_i,s_{i+1})-h(\bar x_i,\bar x_{i+1})\bigr]
=\frac{1}{2\mu}\sum_{i=n}^{m-1}e^{\lb i}(\delta_{i+1}-\delta_i)^2\ge0.
\]
Thus every orbit with $y_0\ne0$ is fixed-endpoint minimizing but lies outside $\wt\cA$. For a related continuous-time rigidity result on the absence of conjugate points, see \cite{ASZ}.
\end{rmk}

\section{variational properties of invariant graphs}\label{s2}

The results of this section apply to the general exact conformally symplectic twist maps introduced in Section~\ref{s1.1}, with $\lambda>0$ fixed.

Following Le Calvez \cite{C,Cal}, we shall use the following notion.

\begin{defn}[Global attractor, Le Calvez \cite{C,Cal}]\label{def:ga}
Let $\Phi$ be a twist map for which some annulus $\mathcal A_M:=\T\times[-M,M]$ is mapped into its interior. Then
\[
  \Om:=\bigcap_{n\ge0}\Phi^{n}(\mathcal A_M)
\]
is independent of $M$; it is called the {\bf global attractor} of $\Phi$. It is compact, connected, invariant and non-contractible; every invariant graph of $\Phi$ is contained in $\Om$, and so is the $\om$-limit set of every orbit contained in $\mathcal A_M$.
\end{defn}

For the CSTM \eqref{eq:cstm} with $\lb>0$ the area is contracted by $e^{-\lb}$. If an invariant graph of \eqref{eq:cstm} exists, it is contained in the global attractor and has particularly fine properties:
\begin{thm}{\label{thm:lip}}
If there exists a graphic invariant circle $\cT$ of \eqref{eq:cstm}, then there exists a $\psi\in {\rm Lip}(\T,\R)$ such that $\cT=\{(x,\psi(x))|x\in\T\}$, with
\[
[\psi]:=\int_{\T}\psi(x) dx=0.
\]
Moreover, for any orbit $\{(x_i,y_i)\}_{i\in\Z}$ in $\cT$, there exists a unique {\bf rotation number}
\beq\label{eq:rot-num}
\rho:=\lim_{i\rightarrow+\infty}\frac{\overline x_i-\overline x_0}{i}.
\eeq
\end{thm}
\proof
By the graphic property, $\psi\in C(\T,\R)$. Since $\cT$ is an invariant graph of a $C^1$ twist map, Birkhoff's graph theorem (see \cite[Chapter 12]{G}) shows that $\psi$ is automatically Lipschitz: a continuous invariant graph of a twist map --- conservative or dissipative alike --- is Lipschitz.

Next, we show that $[\psi]=0$. If not, we suppose $c:=[\psi]$, then $\psi(x)=c+du(x)$ for some $u\in C^{1,1}(\T,\R)$. The invariance of $\cT$ implies for any $x\in\T$ and $n\in\Z$,
\ben
\Phi^{n}\Big(x,c+du(x)\Big)=\Big(x(n), c+du(x(n))\Big)
\een
where $x(0)=x$ and $\{\bar x(n)\}_{n\in\Z}$ denotes the associated critical configuration.
For any $x\in\T$, we have
\[
e^{\lb }(c\overline X+u(X))-u(x)-c\overline x=h(\bar x,\overline X).
\]
Since ${\rm Flux}(\Phi)=0$,
\[
\int_0^1e^{\lb}(c+du(X))dX-(c+du(x)) dx=c(e^\lb-1)=0
\]
which forces $c=0$. Finally, the following diagram commutes:
\[
	\begin{tikzcd}
		\cT \arrow[r, " \Phi"] \arrow[d, "\pi"']
		& \cT\arrow[d, " \pi"] \\
		\mathbb{T} \arrow[r, "\Upsilon"]
		& \mathbb{T}
	\end{tikzcd}
\]
Then $\Upsilon(x):=f(x,du(x))$ becomes a homeomorphism of the circle. Therefore, the rotation number of $\Upsilon$ is uniquely identified as in \eqref{eq:rot-num}.
\qed

\begin{thm}\label{thm:aubry}
$\psi=d u$, where $u(x)$ is the weak KAM solution. Consequently, $\cT$ is globally minimal and coincides with the Aubry set, $\cT=\wt\cA$.
\end{thm}
\begin{proof}
By Theorem~\ref{thm:lip}, $\psi=du$ for a zero-mean $u\in C^{1,1}(\T,\R)$. Set
\[
R(x,X)=h(x,X)+u(x)-e^{\lb}u(X).
\]
In the momentum coordinates $y'=y-du(x)$, $Y'=Y-du(X)$, $R$ generates the
same map and sends $\cT$ to the zero section; indeed,
$-\partial_1R=y'$ and $e^{-\lb}\partial_2R=Y'$. Let
\[
l(x)=\pi\circ\Phi(x,du(x)).
\]
Then $l$ is an orientation-preserving Lipschitz homeomorphism, and
$\partial_1R(x,l(x))=\partial_2R(x,l(x))=0$. Since $x\mapsto R(x,l(x))$ is
absolutely continuous with a.e. derivative
$\partial_1R+\partial_2R\,l'=0$, we have $R(x,l(x))\equiv R_0$.

Since $\partial_{12}R=\partial_{12}h<0$, for $X<l(x)$ and $x_0=l^{-1}(X)$,
\[
\partial_2R(x,X)
=\partial_2R(x_0,X)+\int_{x_0}^{x}\partial_{12}R(s,X)\,ds
=\int_{x_0}^{x}\partial_{12}R(s,X)\,ds<0;
\]
for $X>l(x)$ the same formula gives $\partial_2R(x,X)>0$. Hence
$R(x,\cdot)$ decreases on the left of $l(x)$ and increases on the right, so
\begin{equation}\label{eq:R-min-aubry}
R(x,X)\ge R_0,\qquad\text{with equality iff }X=l(x).
\end{equation}

For an orbit segment $\{\bar x_i\}_{i=n}^m\subset\cT$ and a competitor
$\{s_i\}_{i=n}^m$ with the same endpoints, $h=R-u+e^{\lb}u$ yields
\[
\sum_{i=n}^{m-1}e^{\lb i}\bigl[h(s_i,s_{i+1})-h(\bar x_i,\bar x_{i+1})\bigr]
=\sum_{i=n}^{m-1}e^{\lb i}\bigl[R(s_i,s_{i+1})-R_0\bigr]\ge0,
\]
because the $u$-terms telescope and $\bar x_{i+1}=l(\bar x_i)$. Thus $\cT$ is
globally minimal, and $\psi=du$ is a weak KAM solution.

Finally, $\wt\cA$ is the graph of $du$ over $\cA$, while $\cT$ is its graph
over all of $\T$; hence $\wt\cA\subset\cT$. Since $\cT$ is globally minimal,
$\cT\subset\wt\cA$ by Definition~\ref{def:aubry}, so $\cT=\wt\cA$.
\end{proof}

Theorem~\ref{thm:aubry} shows that an invariant graph is not merely invariant: every orbit on it is globally minimizing. When its rotation number is rational, say $p/q$ in lowest terms, the induced circle homeomorphism $\Upsilon$ of Theorem~\ref{thm:lip} has, by Poincar\'e's theorem, periodic points of period $q$; the corresponding orbits on $\cT$ are $(p,q)$-periodic and, by the same theorem, globally minimizing. The linearization along such an orbit is therefore constrained by the second variation of the action. The next proposition records these constraints.

\begin{prop}\label{prop:cstm-min-periodic}
Let $\Phi$ be an exact CSTM generated by $h$, and let $z_i=(x_i,y_i)$, $i\in\Z$, be a $(p,q)$-periodic orbit, that is, $x_{i+q}=x_i+p$ and $y_{i+q}=y_i$, whose configuration $\{x_i\}$ is fixed-endpoint minimizing, i.e. every finite segment minimizes the weighted action with its endpoints fixed. Put
\[
a_i=\partial_{11}h(x_i,x_{i+1}),\qquad
b_i=\partial_{22}h(x_{i-1},x_i),\qquad
c_i=-\partial_{12}h(x_i,x_{i+1})>0,
\]
\[
D_i=a_i+e^{-\lb}b_i,\qquad
t_i=\frac12\bigl(e^{\lb/2}a_i+e^{-\lb/2}b_i\bigr).
\]
Then $a_i,b_i,c_i,t_i$ are $q$-periodic, and:
\begin{itemize}
\item[(i)] the vertical displacements along the orbit satisfy the Jacobi equation
\[
c_i\,\delta x_{i+1}=D_i\,\delta x_i-e^{-\lb}c_{i-1}\,\delta x_{i-1};
\]
setting $u_i=e^{\lb i/2}\delta x_i$, this becomes
\[
c_i\,u_{i+1}+c_{i-1}\,u_{i-1}=2t_i\,u_i .
\]
\item[(ii)] for every $n\ge1$ the Jacobi matrix
\[
\mathcal H_n=
\begin{pmatrix}
2t_1&-c_1& &\\
-c_1&2t_2&\ddots&\\
&\ddots&\ddots&-c_{n-1}\\
&&-c_{n-1}&2t_n
\end{pmatrix}
\]
is strictly positive definite, not merely positive semi-definite.
\item[(iii)] if $A_j=D\Phi^q(z_j)$, then the two multipliers $s_1,s_2$ of $A_j$ are real and positive, and
\[
s_1s_2=e^{-\lb q},\qquad s_1+s_2\ge 2e^{-\lb q/2},
\]
with equality if and only if $s_1=s_2=e^{-\lb q/2}$. Thus $s_2<1$ corresponds to an attracting node (a node with distinct multipliers if $s_1<s_2$, and an improper (degenerate) node if $s_1=s_2$), $s_2>1$ to a saddle, and $s_2=1$ to the parabolic boundary.
\item[(iv)] for every $j$,
\[
\bigl(D\Phi^q(z_j)\bigr)_{12}>0.
\]
In particular, $D\Phi^q(z_j)$ has no vertical eigendirection.
\end{itemize}
\end{prop}

\begin{proof}
The $q$-periodicity of $a_i,b_i,c_i,t_i$ follows from that of the orbit and of $h$.

(i) The generating relations $y_i=-\partial_1h(x_i,x_{i+1})$ and $y_{i+1}=e^{-\lb}\partial_2h(x_i,x_{i+1})$, linearized along the orbit, give
\[
\delta y_i=-a_i\delta x_i+c_i\delta x_{i+1},
\qquad
\delta y_{i+1}=e^{-\lb}\bigl(-c_i\delta x_i+b_{i+1}\delta x_{i+1}\bigr).
\]
Eliminating $\delta y$ between the second relation and its counterpart at $i+1$ yields the Jacobi equation; the substitution $u_i=e^{\lb i/2}\delta x_i$ gives the second identity.

(ii) Fix $n\ge1$ and consider the segment from $x_0$ to $x_{n+1}$ with its endpoints fixed. Since the orbit is fixed-endpoint minimizing, this segment minimizes the weighted action, so its second variation is nonnegative:
\[
\sum_{i=0}^{n}e^{\lb i}\Bigl[h_{11}\delta x_i^2+2h_{12}\delta x_i\delta x_{i+1}+h_{22}\delta x_{i+1}^2\Bigr]
=e^{-\lb/2}\,\mathbf u^{\top}\mathcal H_n\mathbf u\ \ge 0,
\]
where the $h_{ij}$ are evaluated at $(x_i,x_{i+1})$ and $\mathbf u=(u_1,\dots,u_n)^{\top}$; hence $\mathcal H_n\succeq0$. Let $\Delta_0=1$ and $\Delta_n=\det\mathcal H_n$. The tridiagonal recurrence gives
\[
\Delta_n=2t_n\Delta_{n-1}-c_{n-1}^2\Delta_{n-2}\qquad(n\ge2),
\]
and $\Delta_n\ge0$. If $\Delta_n=0$ for some $n\ge1$, take the smallest such $n$; then $\Delta_{n-1}>0$, and the recurrence at $n+1$ gives $\Delta_{n+1}=2t_{n+1}\Delta_n-c_n^2\Delta_{n-1}=-c_n^2\Delta_{n-1}<0$, contradicting $\mathcal H_{n+1}\succeq0$. Hence $\Delta_n>0$, that is, $\mathcal H_n\succ0$.

(iii) Let $\varphi$ be the solution of the Jacobi equation with $(\varphi_0,\varphi_1)=(0,1)$. Solving the recurrence with the boundary condition $u_0=0$, $u_1=1$ gives
\[
\varphi_{n+1}=\frac{\Delta_n}{c_1c_2\cdots c_n}>0\qquad(n\ge0).
\]
Let $M$ be the monodromy matrix of the recurrence, $M(u_0,u_1)^{\top}=(u_q,u_{q+1})^{\top}$. Its transfer matrices $M_i$ satisfy $\det M_i=c_{i-1}/c_i$, so $\det M=c_0/c_q=1$. Since the coefficients are $q$-periodic, $M^m(u_0,u_1)^{\top}=(u_{mq},u_{mq+1})^{\top}$, and hence
\[
(M^m)_{12}=\varphi_{mq}>0\qquad(m\ge1).
\]
Write $\tau=\operatorname{tr}M$. If $\tau\le0$, Cayley--Hamilton and $\det M=1$ give $(M^2)_{12}=\tau M_{12}\le0$, a contradiction. If $0<\tau<2$, write $\tau=2\cos\theta$ with $0<\theta<\pi/2$; again by Cayley--Hamilton, $(M^m)_{12}=\sin(m\theta)M_{12}/\sin\theta$, which is nonpositive for some positive integer $m$, again a contradiction. Hence $\tau\ge2$.

It remains to relate $M$ to $A_j=D\Phi^q(z_j)$. In the coordinates $(\delta x_i,\delta x_{i+1})$ the linearized recurrence acts by a matrix $N$, which is similar to $A_j$: the change of variables $\delta y_i=-a_i\delta x_i+c_i\delta x_{i+1}$ is invertible because $c_i>0$. Since $u_i=e^{\lb i/2}\delta x_i$ and the weights are $q$-periodic, $M=e^{\lb q/2}DND^{-1}$ with $D=\operatorname{diag}(1,e^{\lb/2})$, so $M$ is similar to $e^{\lb q/2}N$. Therefore
\[
s_1+s_2=e^{-\lb q/2}\operatorname{tr}M\ \ge\ 2e^{-\lb q/2},\qquad s_1s_2=\det A_j=e^{-\lb q}.
\]
The discriminant $(s_1+s_2)^2-4e^{-\lb q}\ge0$ shows that $s_1,s_2$ are real, and they are positive since $s_1+s_2>0$ and $s_1s_2>0$; the inequality is an equality if and only if $s_1=s_2=e^{-\lb q/2}$. The classification follows from $0<s_1\le s_2$ and $s_1s_2=e^{-\lb q}$.

(iv) Take the vertical initial displacement $\delta x_0=0$, $\delta y_0=1$. Then $\delta x_1=1/c_0$, so $(u_0,u_1)=(0,e^{\lb/2}/c_0)$, and by linearity of the Jacobi equation
\[
u_q=\frac{e^{\lb/2}}{c_0}\,\varphi_q,
\qquad\text{hence}\qquad
\bigl(D\Phi^q(z_0)\bigr)_{12}=\delta x_q=e^{-\lb(q-1)/2}\frac{\varphi_q}{c_0}>0.
\]
The same argument applies to every cyclic shift of the orbit. If $(0,1)^{\top}$ were a vertical eigenvector of $D\Phi^q(z_j)$, its first component after applying $D\Phi^q(z_j)$ would vanish, contradicting the strict positivity just proved. Hence no eigendirection is vertical.
\end{proof}

\begin{cor}\label{cor:graph-periodic-linear}
Let $\cT$ be an invariant graph whose rotation number is rational, say $p/q$ in lowest terms. Then every $(p,q)$-periodic point $z_0\in\cT$ satisfies the conclusions of Proposition~\ref{prop:cstm-min-periodic}. In particular, $D\Phi^q(z_0)$ has positive real multipliers and no vertical eigendirection, even if $\cT$ is not differentiable at $z_0$.
\end{cor}

\begin{proof}
By Theorem~\ref{thm:aubry}, $\cT=\wt\cA$ and every orbit on $\cT$ is globally minimal, hence fixed-endpoint minimizing (Lemma~\ref{lem:aubry-fixed-endpoint}). Applying Proposition~\ref{prop:cstm-min-periodic} to the periodic orbit through $z_0$ gives the claim.
\end{proof}

\section{The breakdown of invariant graphs for dissipative standard maps}\label{s3}

Let $\mu=(1-e^{-\lb})/\lb$ and suppose $V\in C^2(\T,\R)$, so that $V'\in C^1(\T,\R)$ is periodic with zero mean. Throughout this section, no special shape of $V$ or $V'$ is imposed unless explicitly stated. The upper-bound result requires $V$ to be nonconstant (equivalently, $V'\not\equiv0$), whereas the lower-bound estimates also cover general periodic potentials. We fix $\lambda>0$ throughout. For each chosen $\omega$, we vary $\epsilon\ge0$ to study the existence and breakdown of invariant graphs. We also compare the resulting thresholds across different values of $\omega$, always keeping $\lambda$ fixed.

The {\bf dissipative standard map (DSM)} is given by
\be\label{eq:stan}
\Phi_{\lb,\om,\eps}:\left\{
\begin{aligned}
&X=x+\mu y+\om(\in \R),\\
&Y=e^{-\lb}y+\eps e^{-\lb}V'(X)
\end{aligned}
\right.
\ee
with the generating function
\be
h(x,X)=\frac{1}{2\mu}|X-x-\om|^2+\eps V(X),\quad\eps \ge 0.
\ee

This system possesses a global attractor. Indeed, since $V'$ is bounded, \eqref{eq:stan} gives
\[
|y_n|\le e^{-\lb n}|y_0|+\eps\|V'\|_\infty\frac{e^{-\lb}}{1-e^{-\lb}}\qquad(n\ge0),
\]
so every orbit is contained in some annulus $\mathcal A_M$ with $M\ge\eps\|V'\|_\infty e^{-\lb}/(1-e^{-\lb})$, while $\Phi(\mathcal A_M)\subset\operatorname{int}\mathcal A_M$ for such $M$. Hence, by Definition~\ref{def:ga}, the global attractor exists and {\it contains the $\om$-limit set of every orbit}. It also contains an invariant circle.

For sufficiently small values of the parameter $\epsilon$, the invariant circle exists and is graphic. We will prove this in \ref{ss:lb}.

In this section we study the breakdown threshold
\[
\eps_*(\lb,\om):=\inf\{\eps\ge0:\ \Phi_{\lb,\om,\eps}\ \text{has no invariant graph}\},
\]
that is, the infimum of the values of $\eps$ for which the invariant graph associated with $(\lb,\om)$ does not exist.

\subsection{Herman's non-existence of graphs} In this subsection we assume that $V\in C^2(\T,\R)$ is nonconstant, or equivalently that $V'\in C^1(\T,\R)$ is nonzero, periodic, and has zero mean. The estimates depend on the extremal values of $V''$, but not on a particular trigonometric form. The first result is an upper bound for $\eps_*(\lb,\om)$.

Here $\lambda>0$ is fixed, $\omega\in\R$ is arbitrary, and $\epsilon$ varies. The bound $K_H(\lambda)$ is independent of $\omega$, so the non-existence criterion holds uniformly over all drift parameters $\omega$.

\begin{thm}\label{thm:K_H}
    For any nonconstant $V\in C^2(\T,\R)$ with $\lb>0$ fixed, we can always find a $0<K_H(\lb)<+\infty$ depending on $V''$ and $ \lb$, such that for any $\eps\geq K_H(\lb)$, the associated \eqref{eq:stan} admits no invariant graphs.
\end{thm}

\proof Suppose the graph of $\psi:\T\rightarrow\R$ is invariant, i.e.\ $\cT:=\{(x,\psi(x))\mid x\in\T\}$ is an invariant circle of $\Phi_{\lb,\om,\eps}$; then
\be\label{eq:graph}
    \psi(x+\mu\psi(x)+\om) = e^{-\lb}\Big(\psi(x)+\epsilon V'(x+\mu\psi(x)+\om)\Big).
\ee
By Theorem~\ref{thm:lip}, $\psi$ is a Lipschitz continuous function.

Because of the orientation-preserving property of \eqref{eq:stan} (restricted to the invariant graph), the function $f:=x+\mu\psi(x)+\om$ is invertible when lifted to $\R$. Note that we do not presuppose $\psi$ to be differentiable in advance. Since $\psi$ is necessarily Lipschitz continuous, Rademacher's theorem implies that $\psi$ is differentiable almost everywhere, and so is $f$. Since $V\in C^2(\T,\R)$, the function $\eps V'$ is of class $C^1$; consequently all the terms in \eqref{eq:graph} are differentiable almost everywhere, and the chain rule may be applied almost everywhere.

We can rewrite \eqref{eq:graph} as
\[
\psi\circ f=e^{-\lb}\psi\circ f^{-1}\circ f+e^{-\lb}\epsilon V'\circ f,\quad\forall x\in[0,1)
\]
and then further as
\be
e^\lb f+f^{-1}=(e^\lb+1)id+(e^\lb-1)\om+\mu \epsilon V'.
\ee
For convenience we reduce it to
\[
e^{\lb/2}f+e^{-\lb/2}f^{-1} = 2\cosh\Big(\frac \lb2\Big)id+ 2\sinh\Big(\frac \lb2\Big)\cdot\Big(\om+\frac {e^{-\lb}}\lb \epsilon V'\Big).
\]
Differentiating both sides of this equality, we obtain
\be\label{eq:sh-gra}
\frac 12 \Big(e^{\lb/2}f'+\frac1{e^{\lb/2}f'\circ f^{-1}}\Big) = \cosh\Big(\frac \lb2\Big)+\sinh\Big(\frac \lb2\Big)\frac {e^{-\lb}}\lb \epsilon V'', \quad \text{a.e. }x\in\T.
\ee
Suppose
\[
M:=\max\big\{\ess\sup_{\overline x\in\R}e^{\lb/2}f',\ess\sup_{\overline x\in\R} \frac 1{e^{\lb/2}f'} \big\},\quad m:=\min\big\{\ess\inf_{\overline x\in\R}e^{\lb/2}f',\ess\inf_{\overline x\in\R} \frac1{e^{\lb/2}f'} \big\},
\]
then $m\geq1/M$. Let
     \[
     v_- := \ess\inf_{\overline{x}\in\R} V''(\overline{x}),\quad v_+ := \ess\sup_{\overline{x}\in\R} V''(\overline{x}).
     \]
Since $V'$ is periodic, $v_-\le0$, and hence

\be\label{eq:M1}
\cosh\Big(\frac \lb2\Big)+\epsilon\sinh\Big(\frac \lb2\Big)\frac {e^{-\lb}}\lb v_-\geq \frac 1M
\ee
and
\be\label{eq:M2}
\frac 12 \Big(M+\frac1M\Big) \leq \cosh\Big(\frac \lb2\Big)+\eps\sinh\Big(\frac \lb2\Big)\frac {e^{-\lb}}\lb v_+
\ee
due to (\ref{eq:sh-gra}). Combining inequalities (\ref{eq:M1}) and (\ref{eq:M2}), we finally obtain a criterion
\be\label{eq:gra-cri}
\left\{
\begin{aligned}
M\geq & \Big(\cosh\Big(\dfrac \lb2\Big)+\eps\sinh\Big(\dfrac \lb2\Big)\dfrac{v_-}{\lb e^{\lb}}\Big)^{-1},\\
M\geq & \cosh\Big(\dfrac \lb2\Big)+\eps\sinh\Big(\dfrac \lb2\Big)\dfrac{v_+}{\lb e^{\lb}} -\sqrt{\Big[\cosh\Big(\dfrac \lb2\Big)+\eps\sinh\Big(\dfrac \lb2\Big)\dfrac{v_+}{\lb e^{\lb}}\Big]^2-1},\\
M\leq & \cosh\Big(\dfrac \lb2\Big)+\eps\sinh\Big(\dfrac \lb2\Big)\dfrac{v_+}{\lb e^{\lb}} +\sqrt{\Big[\cosh\Big(\dfrac \lb2\Big)+\eps\sinh\Big(\dfrac \lb2\Big)\dfrac{v_+}{\lb e^{\lb}}\Big]^2-1}.
\end{aligned}
\right.
\ee
This formula is equivalent to
\be\label{eq:final-est}
\left\{
\begin{aligned}
(R1) \quad & 0 \le -\eps v_- < \lb e^{\lb}\coth\Big(\frac{\lb}{2}\Big) \\
(R2) \quad & \eps\sinh\Big(\dfrac \lb2\Big)\dfrac{(v_+-v_-)}{\lb e^{\lb}} \le \sqrt{\Big[\cosh\Big(\dfrac \lb2\Big)+\eps\sinh\Big(\dfrac \lb2\Big)\dfrac{v_+}{\lb e^{\lb}}\Big]^2-1}
\end{aligned}
\right.
\ee

In summary, $\eps$ must satisfy the following condition:
\begin{equation}\label{eq:K_H_V}
    \eps \le \lb e^{\lb}\frac{\sqrt{v_-^2-2v_-v_++v_+^2 \coth^2(\frac{\lambda}{2})}+v_+\coth(\frac{\lambda}{2})}{v_-^2-2v_-v_+}=:K_H(\lb)
\end{equation}
The right-hand side remains finite for $\lb>0$ and every nontrivial $V(x)$.
\qed

\begin{cor}\label{cor:herman}
Suppose $V'(x)=\sin(2\pi x)$. Then for all $\eps\geq K_H(\lb)$ with
\be\label{eq:final-est-sin}
K_H(\lb):=\frac{\lb e^\lb}{2\pi\cdot3}\Big\{\coth\Big(\frac \lb2 \Big)+\sqrt{\coth^2\Big(\frac \lb2 \Big)+3}\Big\}
\ee
the DSM \eqref{eq:stan} admits no invariant graphs. Besides,
\[ \begin{split}
\lim_{\lb \to 0^+} K_H(\lb) = & \frac {1}{2\pi} \cdot \frac 23 \lb e^\lb\coth\Big(\frac \lb2\Big)=\frac {1}{2\pi} \cdot \frac 43, \\ K_H(\lb)\sim & \frac{\lb e^\lb}{2\pi}\coth\frac{\lb}{2}\sim \frac{\lb e^\lb}{2\pi}\qquad (\lb\to+\infty).
\end{split} \]
These estimates coincide with the computations in \cite{Mat2} for area-preserving standard maps.
\end{cor}

\subsection{Lower-bound estimate for the non-existence of graphs}\label{ss:lb}

Here $V\in C^2(\T,\R)$ remains arbitrary, including potentials with several critical points; $V'\in C^1$ may have any periodic zero-mean profile. No nondegeneracy or sine-profile assumption is needed for the iteration and the general estimates below. Only the explicit identity for $K_0(\lambda,0)$ is specialized to $V'(x)=\sin(2\pi x)$.

Fix $\lambda>0$ and $\omega\in[0,1/2]$, and vary $\epsilon\ge0$; we abbreviate $\Phi_{\lambda,\omega,\epsilon}$ by $\Phi_\epsilon$. When comparing the thresholds $K_f(\lambda,\omega)$ for different values of $\omega$, we always use the same fixed $\lambda$.

Instead of iterating only the zero section $\Gamma_0=\{(x,0)\in T^*\T\mid x\in\T\}$, we start from an arbitrary graph. For $f\in C(\T,\R)$ put
\[
\Gamma_{f,0}=\{(x,f(x)):x\in\T\},\qquad \Gamma_{f,n+1}=\Phi_\eps(\Gamma_{f,n})\quad(n\in\N),
\]
and write $\Gamma_n=\Gamma_{0,n}$ when $f\equiv0$; whenever $\Gamma_{f,n}$ is a graph, $\psi_{f,n}$ denotes the function describing it.

\begin{defn}\label{def:K_f}
For $f\in C(\T,\R)$ set
\beq\label{eq:low-b}
K_f(\lambda,\omega) := \inf\bigl\{ \epsilon\ge 0 \mid \exists\, n\in\mathbb{N} \text{ such that } \Gamma_{f,n} \text{ is not a graph} \bigr\},
\eeq
with the convention $\inf\varnothing = +\infty$; for $f\equiv0$ we write $K_0=K_f$. By the definition of the infimum, for every $0\le \epsilon<K_f$ all iterates $\Gamma_{f,n}$ are graphs. In particular, $K_f=0$ if the first iterate $\Gamma_{f,1}$ is already not a graph at $\eps=0$.
\end{defn}

Notice that $K_f(\lb,\om)$ can be $+\infty$, e.g., when $V$ is constant. Therefore we only consider nontrivial $V(x)$, for which $K_f(\lb,\om)\le\eps_*(\lb,\om)\le K_H(\lb)<+\infty$ by Proposition~\ref{prop:Kf-bound} below.

\begin{lem}[Uniform Lipschitz bound]\label{lem:ulb}
If $\epsilon<K_f$, then for every $n\ge1$ the function $\psi_{f,n}$ describing $\Gamma_{f,n}$ is Lipschitz continuous and satisfies
\[
\operatorname{Lip}(\psi_{f,n})\le L_\epsilon:=\max\Bigl\{\frac1\mu,\;\frac{e^{-\lambda}}{\mu}+\epsilon e^{-\lambda}\|V''\|_\infty\Bigr\}.
\]
\end{lem}

\begin{proof}
We drop the subscript $f$ and write $\Gamma_n=\Gamma_{f,n}$, $\psi_n=\psi_{f,n}$. Because $\Gamma_n$ is a graph, the induced circle map
\(T_n(x)=x+\mu\psi_n(x)+\omega\) must be non-decreasing; its derivative satisfies $T_n'(x)=1+\mu\psi_n'(x)\ge 0$ almost everywhere, hence $\psi_n'\ge -1/\mu$ a.e.
From~\eqref{eq:stan} we obtain, for almost every $X$,
\[
\psi_{n+1}(X)=e^{-\lambda}\psi_n(T_n^{-1}(X))+\epsilon e^{-\lambda}V'(X).
\]
All involved maps are Lipschitz; therefore the chain rule applies almost everywhere and yields
\[
\psi_{n+1}'(X)=e^{-\lambda}\frac{\psi_n'(x)}{1+\mu\psi_n'(x)}+\epsilon e^{-\lambda}V''(X),\qquad x=T_n^{-1}(X).
\]
Since $1+\mu\psi_n'(x)\ge0$ a.e., for every $t\ge-1/\mu$ one has $\dfrac{t}{1+\mu t}\le\dfrac1\mu$, whence
\[
\psi_{n+1}'(X)\le\frac{e^{-\lb}}{\mu}+\eps e^{-\lb}\|V''\|_\infty \qquad\text{a.e.}
\]
For the lower bound we use the hypothesis $\eps<K_f$: then $\Gamma_{n+1}$ is a graph, so the induced circle map $T_{n+1}(x)=x+\mu\psi_{n+1}(x)+\om$ is non-decreasing, and therefore $\psi_{n+1}'(X)\ge-1/\mu$ a.e. Combining the two bounds yields the claimed Lipschitz constant.
\end{proof}

\begin{lem}[Uniform boundedness]\label{lem:ub}
Under the same hypotheses, for every fixed $f\in C(\T,\R)$ there is a constant $M_\epsilon(f)$, independent of $n$, such that $\|\psi_{f,n}\|_\infty\le M_\epsilon(f)$.
\end{lem}

\begin{proof}
In the notation of the preceding proof, from the expression of $\psi_{n+1}$ we have

 \ben
 \|\psi_{n+1}\|_\infty & \le & e^{-\lambda}\|\psi_n\|_\infty+\epsilon e^{-\lambda}\|V'\|_\infty\\
 & \le & e^{-(n+1)\lb}\|\psi_0\|_\infty+\epsilon e^{-\lambda}\frac{1-e^{-(n+1)\lb}}{1-e^{-\lb}}\|V'\|_\infty\\
 & \le & \|\psi_0\|_\infty+\frac{\epsilon e^{-\lambda}}{1-e^{-\lb}}\|V'\|_\infty.
 \een
so the sequence is bounded uniformly in $n$ by $M_\epsilon(f):=\|f\|_\infty+\epsilon e^{-\lambda}\|V'\|_\infty/(1-e^{-\lambda})$, which depends on $f$ but not on $n$.
\end{proof}

\begin{thm}[Existence and uniqueness of the invariant graph]\label{thm:e_u_i}
For every $f\in C(\T,\R)$ and every $\epsilon\in[0,K_f)$, the map $\Phi_\epsilon$ possesses a unique invariant circle, and this circle is a Lipschitz graph.
\end{thm}

\begin{proof}
Fix $f$ and use the notation of Lemma~\ref{lem:ulb}. Lemmas~\ref{lem:ulb} and~\ref{lem:ub} show that, for fixed $f$, $\{\psi_n\}$ is uniformly bounded in $n$ and uniformly Lipschitz, hence pre-compact in $C^0(\mathbb{T})$ by Arzel\`a--Ascoli.

Consider the area between two successive graphs,
\[
A_n=\int_0^1|\psi_{n+1}(x)-\psi_n(x)|\,\mathrm{d}x.
\]
The set $S_n=\{(x,y)\in T^*\T\mid \text{$y$ lies between $\psi_n(x)$ and $\psi_{n+1}(x)$}\}$ has Lebesgue measure $A_n$. Because $\Phi_\epsilon$ is a bi-Lipschitz homeomorphism with constant Jacobian determinant $\det D\Phi_\epsilon = e^{-\lambda}$ almost everywhere, it maps $S_n$ onto $S_{n+1}$ and the area formula gives
\[
A_{n+1}=e^{-\lambda}A_n,\qquad\text{hence}\qquad A_n=e^{-\lambda n}A_0\to0.
\]
In particular, $\|\psi_{n+1}-\psi_n\|_{L^1}\to0$.

If $\{\psi_n\}$ were not Cauchy in $C^0$, we could extract two subsequences converging uniformly to distinct limits $\psi$ and $\tilde\psi$ with $\|\psi-\tilde\psi\|_\infty\ge\delta>0$. However, the $L^1$ distance between these limits must be zero because the successive $L^1$ differences tend to zero, forcing $\psi=\tilde\psi$ by continuity. Thus $\{\psi_n\}$ is Cauchy and converges uniformly to some $\psi\in C^0(\mathbb{T})$. The uniform Lipschitz bound passes to the limit, so $\psi$ is Lipschitz.

To check invariance, pass to the limit in the functional equation
\[
\psi_{n+1}(X)=e^{-\lambda}\psi_n(T_n^{-1}(X))+\epsilon e^{-\lambda}V'(X)
\]
using the uniform convergence of $\psi_n$ and of the inverses $T_n^{-1}$ (which follows from the uniform bi-Lipschitz estimates). We obtain $\psi=e^{-\lambda}\psi\circ T^{-1}+\epsilon e^{-\lambda}V'$, i.e.\ $\Phi_\epsilon(\Gamma)=\Gamma$ for $\Gamma=\{(x,\psi(x))\mid x\in\T\}$.

Uniqueness. We first recall the argument for two invariant graphs. If $\psi$ and $\tilde\psi$ are two invariant graphs, the area $A=\int_0^1|\psi-\tilde\psi|dx$ between them satisfies $A=e^{-\lb}A$, hence $A=0$ and $\psi=\tilde\psi$; so there is at most one invariant graph.

Next we show that two invariant circles cannot intersect. Suppose that $C_1\ne C_2$ are invariant circles with $C_1\cap C_2\ne\varnothing$, and put $K=C_1\cup C_2$. Since $\Phi(K)=K$, $\Phi$ permutes the connected components of $T^*\T\setminus K$; since $\Phi$ and $\Phi^{-1}$ map bounded sets onto bounded sets, it permutes the bounded components among themselves. Let $S$ be the union of the bounded components of $T^*\T\setminus K$. Because the two circles are distinct, $S$ has positive area; and $\Phi(S)=S$, so $\operatorname{area}(S)=e^{-\lb}\operatorname{area}(S)=0$, a contradiction. Hence two distinct invariant circles are disjoint.

Finally, if $C_1\ne C_2$ are disjoint invariant circles, the annulus between them is invariant, so its area also satisfies $\operatorname{area}=e^{-\lb}\operatorname{area}$, hence is zero, again a contradiction. Therefore there is at most one invariant circle; since $\Gamma$ is one, it is the unique invariant circle.
\end{proof}

\begin{prop}\label{prop:Kf-bound}
For every $f\in C(\T,\R)$ one has $K_f(\lb,\om)\le\eps_*(\lb,\om)$, and therefore
\begin{equation}\label{eq:supKf}
\sup_{f\in C(\T,\R)}K_f(\lb,\om)\ \le\ \eps_*(\lb,\om).
\end{equation}
The supremum on the left is a lower bound for $\eps_*(\lb,\om)$ obtained from the iterates of a single initial graph, and it need not coincide with $\eps_*(\lb,\om)$.
\end{prop}

\begin{proof}
Let $\eps<K_f$. Then all $\Gamma_{f,n}$ are graphs, and Theorem~\ref{thm:e_u_i} provides an invariant graph for $\Phi_\eps$; hence $\eps\le\eps_*$. Since this holds for every $\eps<K_f$, we obtain $K_f\le\eps_*$, and \eqref{eq:supKf} follows by taking the supremum over $f$.
\end{proof}

\begin{prop}\label{prop:plb}
Let $f\equiv0$ be the zero section, and assume $\|V''\|_\infty>0$. Then
\[
K_0(\lb,\om) \ge \frac{\lambda e^{\lambda}}{\|V''\|_\infty}\tanh{\frac{\lb}{4}}.
\]
The bound is specific to the zero section: the induction below starts from $\psi_0\equiv0$, and for a general initial graph the corresponding base case does not necessarily exist (indeed $K_f=0$ whenever the first iterate is already not a graph).
\end{prop}

\begin{proof}
Set $\delta=e^{-\lambda/2}\in(0,1)$. We prove that whenever
\[
\epsilon\le K_*:=\frac{\lambda e^{\lambda}}{\|V''\|_\infty}\tanh{\frac{\lb}{4}},
\]
the iterates satisfy the stronger inequality $1+\mu\psi_n'\ge \delta$ almost everywhere for all $n$; in particular each $\Gamma_n$ is a graph. Consequently such $\epsilon$ do not belong to the set whose infimum defines $K_0$, and therefore $K_0\ge K_*$. The estimate uses the zero section through the base case $\psi_0\equiv0$; for an arbitrary $f$ nothing is asserted about $K_f$.

We proceed by induction on $n$. For $n=0$, $\psi_0\equiv0$ gives $1+\mu\psi_0'=1\ge\delta$.
Assume $1+\mu\psi_n'\ge\delta$ a.e.\ for some $n$. Using the derivative formula from the proof of Lemma~\ref{lem:ulb},
\[
1+\mu\psi_{n+1}'(X)=1+e^{-\lambda}-\frac{e^{-\lambda}}{1+\mu\psi_n'(x)}+\mu\epsilon e^{-\lambda}V''(X),\qquad x=T_n^{-1}(X).
\]
From the induction hypothesis, $-\frac{1}{1+\mu\psi_n'}\ge -\frac1\delta$. Together with $V''(X)\ge -\|V''\|_\infty$ we obtain, for almost every $X$,
\[
1+\mu\psi_{n+1}'(X)\ge 1+e^{-\lambda}-\frac{e^{-\lambda}}{\delta}-\mu\epsilon e^{-\lambda}\|V''\|_\infty.
\]
The right-hand side is a decreasing function of $\epsilon$. Evaluating at $\epsilon=K_*$ gives exactly $\delta$. Hence, for any $\epsilon\le K_*$,
\[
1+\mu\psi_{n+1}'(X)\ge \delta\quad\text{a.e.},
\]
completing the induction. The inequality is not always strict: Corollary~\ref{cor:K0zero} gives an example in which equality holds.
\end{proof}

\begin{cor}\label{cor:K0zero}
For $V'(x)=\sin(2\pi x)$ and the zero section $f\equiv0$, $K_0(\lambda,0)=K_*=\frac{\lambda e^{\lambda}}{2\pi}\tanh{\frac{\lb}{4}}$. Hence
\[
K_0(\lambda,0)\ \le\ K_0(\lambda,\omega)\ \le\ \eps_*(\lambda,\omega)\ \le\ K_H(\lambda),
\]
where the first two inequalities concern the zero section (for another initial graph only $K_f\le\eps_*$ is available, and $K_f$ can be smaller than $K_*$), while the last one is the general bound of Section~\ref{s3}.
\end{cor}

\begin{proof}
For this potential and $f\equiv0$, the recursion preserves the symmetry $\psi_n(1-x)=-\psi_n(x)$, and one checks inductively that $\psi_n'$ attains its maximum at $X=0$ and its minimum at $X=1/2$. Put $g_n=1+\mu\psi_n'(1/2)$. Since $T_n(1/2)=1/2$ and $V''(1/2)=-2\pi$, the derivative formula gives
\[
g_{n+1}=1+e^{-\lambda}-\frac{e^{-\lambda}}{g_n}-2\pi\mu\epsilon e^{-\lambda}.
\]
At $\epsilon=K_*$, $g\equiv\delta=e^{-\lambda/2}$ is the unique fixed point, so $g_n\downarrow\delta>0$ and all $\Gamma_n$ are graphs. For $\epsilon>K_*$ the recursion has no positive fixed point, so $g_n$ eventually becomes negative and some $\Gamma_n$ is not a graph. Hence $K_0(\lambda,0)=K_*$.
\end{proof}

\begin{rmk}\label{rmk:KMW}
By Section~\ref{subsec:bifurcation}, for $V'(x)=\sin(2\pi x)$ and $\omega=0$ one can also prove
\[
K_*=K_0(\lambda,0)=\eps_*(\lambda,0).
\]
\end{rmk}

In the following construction, we continue to fix $\lambda>0$, but unlike the preceding estimates we allow the potential $V(\cdot,\epsilon)$ itself to vary with $\epsilon$. The resulting non-existence statement holds for arbitrary $\omega$.

\begin{cor}
There exists a parametrized function $V(\cdot,\eps)$ whose $C^{2-\dt}$ norm (for any $0<\dt<1$) tends to $0$ as $\eps\to0$, such that the associated (\ref{eq:stan}) has no invariant graph.
\end{cor}

\begin{proof}
    Fix any $\delta\in(0,1)$. Choose $b>1/\delta$, and then choose $a$ such that
    \[
    b(2-\delta)<a<2b-1;
    \]
    such an $a$ exists, since $b(2-\delta)<2b-1$ is equivalent to $b\delta>1$.

    Let $N_{\epsilon}=\lfloor\epsilon^{-b}\rfloor$ and define
    \[
    V_{\epsilon}(x)=\epsilon^{a}\sin\bigl(2\pi N_{\epsilon}x\bigr).
    \]
    Then
    \[
    \|V_{\epsilon}\|_{\infty}\sim\epsilon^{a},\qquad\|V_{\epsilon}'\|_{\infty}\sim\epsilon^{a}N_{\epsilon}\sim\epsilon^{a-b}.
    \]

    We next estimate the H\"older seminorm $[V_{\epsilon}']_{1-\delta}$. Writing
    \[
    V_\eps'(x)=\underbrace{2\pi\eps^a N_\eps}_{G_\eps}\cdot\underbrace{\cos(2\pi N_\eps x)}_{g(x)}=:G_\eps\cdot g(x),
    \]
    we have $[V_{\epsilon}']_{1-\delta}=G_{\epsilon}\,[g]_{1-\delta}$, where
    \[
    [g]_{1-\delta}=\sup_{x\neq y}\frac{|g(x)-g(y)|}{\operatorname{dist}(x,y)^{1-\delta}}.
    \]
    If $\operatorname{dist}(x,y)\geq\frac{1}{2N_{\epsilon}}$, then
    \[
    \frac{|g(x)-g(y)|}{\operatorname{dist}(x,y)^{1-\delta}}\leq 2\bigl(2N_{\epsilon}\bigr)^{1-\delta}=2^{2-\delta}N_{\epsilon}^{1-\delta}.
    \]
    If $\operatorname{dist}(x,y)\leq\frac{1}{2N_{\epsilon}}$, the mean value
    theorem together with $|g'|\leq2\pi N_{\epsilon}$ gives
    $|g(x)-g(y)|\leq2\pi N_{\epsilon}\operatorname{dist}(x,y)$, and hence
    \[
    \frac{|g(x)-g(y)|}{\operatorname{dist}(x,y)^{1-\delta}}\leq2\pi N_{\epsilon}\operatorname{dist}(x,y)^{\delta}\leq2^{1-\delta}\pi N_{\epsilon}^{1-\delta}.
    \]
  Combining the two cases, we obtain
    $[g]_{1-\delta}\leq2^{1-\delta}\pi N_{\epsilon}^{1-\delta}$, and therefore
    \[
    [V_{\epsilon}']_{1-\delta}\leq2^{2-\delta}\pi^{2}\epsilon^{a}N_{\epsilon}^{2-\delta}=C\,\epsilon^{a-b(2-\delta)}
    \]
    for some constant $C$ independent of $\epsilon$. Since $a>b(2-\delta)$, it
    follows that
    \[
    \|V_{\epsilon}\|_{C^{2-\delta}}=\|V_{\epsilon}\|_{\infty}+\|V_{\epsilon}'\|_{\infty}+[V_{\epsilon}']_{1-\delta}\longrightarrow0\quad\text{as }\epsilon\to0.
    \]
    On the other hand,
    $\|V_{\epsilon}''\|_{\infty}\sim\epsilon^{a}N_{\epsilon}^{2}\sim\epsilon^{a-2b}$, so by \eqref{eq:K_H_V} the Herman threshold satisfies
    \[
    K_{H}(\epsilon)=C'(\lambda)\,\epsilon^{2b-a},\qquad C'(\lambda)>0.
    \]
    Since $2b-a>1$, we have
    \[
    \epsilon\geq K_{H}(\epsilon)\iff C'(\lambda)\,\epsilon^{2b-a-1}\leq1,
    \]
    and the left-hand side tends to $0$ as $\epsilon\to0$; hence $\epsilon\geq K_{H}(\epsilon)$ holds for all sufficiently small $\epsilon$. Theorem~\ref{thm:K_H} then implies that the associated map \eqref{eq:stan} admits no invariant graph.

    Finally, the exponent $2-\delta$ is sharp. Indeed, since
    \[
    \|V_{\epsilon}''\|_{\infty}\leq\|V_{\epsilon}\|_{C^{2}}\leq\|V_{\epsilon}\|_{C^{2+\delta}},
    \]
    if either $\|V_{\epsilon}\|_{C^{2}}\to0$ or
    $\|V_{\epsilon}\|_{C^{2+\delta}}\to0$ as $\epsilon\to0$, then
    $\|V_{\epsilon}''\|_{\infty}\to0$, and consequently
    \[
    K_{0}\geq\frac{\lambda e^{\lambda}}{\|V''_\eps\|_{\infty}}\tanh\frac{\lambda}{4}\longrightarrow+\infty,
    \]
    which means that the invariant graph persists. Hence the $C^{2-\delta}$-norm cannot be replaced by the $C^{2}$- or the $C^{2+\delta}$-norm.
\end{proof}

\begin{rmk}\label{rmk:SW}
The exponent $2-\dt$ in the preceding corollary matches \cite[Theorem 1]{SW}: the perturbation in
\eqref{eq:stan} is $\eps e^{-\lb}V'$, so $V\in C^{2-\dt}$ amounts to a perturbation in
$C^{1-\dt}$, exactly the class in which \cite{SW} destroys all $C^0$-invariant graphs, which is sharp by
the NHIM theorem.
\end{rmk}

\section{Mechanism of breakdown \& the drift of rotation numbers}\label{s4}

Throughout this section $\lambda>0$ is fixed; we vary $\eps$ for a chosen $\om$ and compare thresholds and rotation numbers as $\om$ varies. Numerical experiments use $\lambda=2$ except in Section~\ref{subsec:sga}. We assume the generic case $V\in C^2$ with $V'$ having a unique maximum and minimum; no nondegeneracy of these extrema is required. This concerns $V'$, not the extrema of $V$. Results not using uniqueness are stated for general $V$; the smooth-linearization argument needs extra smoothness and non-resonance, while the numerics use $V'(x)=\sin(2\pi x)$.

\subsection{Numerical threshold of the breakdown}

\begin{figure}[htbp]
\begin{center}
\includegraphics[width=13cm]{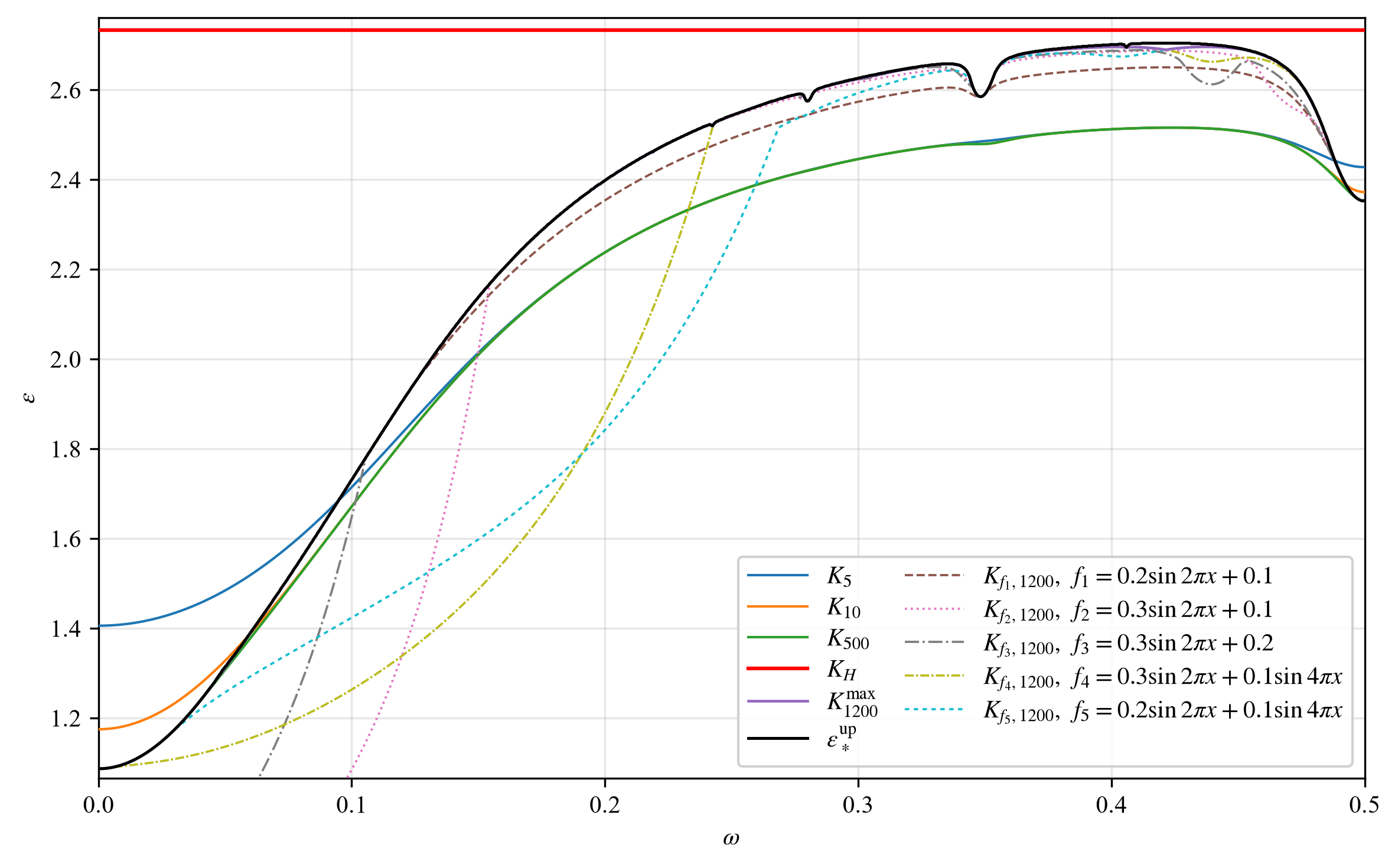}
\caption{With $V'(x)=\sin(2\pi x)$ and $\lambda=2$: the solid curves are the zero-section thresholds $K_5$, $K_{10}$, $K_{500}$ and the red line is the constant bound $K_H$ of Theorem~\ref{thm:K_H}; the purple curve is $K_{1200}^{\text{max}}=\max_iK_{f_i,1200}$ and the dashed curves are the five thresholds $K_{f_i,1200}$ of Appendix~\ref{app:error}; the black curve is $\eps_*^{\text{up}}$. The threshold $\eps_*(\om)$ lies in the band between $K_{1200}^{\text{max}}$ and $\eps_*^{\text{up}}$.}
\label{fig:omega-epsilon-bound}
\end{center}
\end{figure}

In this subsection we set $\lambda=2$. For each $\omega\in[0,1/2]$ we vary $\epsilon$ to estimate the breakdown threshold and its finite-iteration approximations. Figure~\ref{fig:omega-epsilon-bound} then compares these estimates as $\omega$ varies, with $\lambda$ kept fixed throughout.

The computations below use $V'(x)=\sin(2\pi x)$ and $\lb=2$. By Theorem~\ref{thm:K_H}, $\eps_*(\lb,\om)\le K_H$ for every $\om$ (the horizontal line in Fig.~\ref{fig:omega-epsilon-bound}), while Proposition~\ref{prop:Kf-bound} supplies the lower bounds $K_f\le\eps_*$. Since $K_f$ is not directly computable, we use the finite-iteration threshold
\[
K_{f,M}(\lb,\om):=\inf\bigl\{\eps\ge0:\ \Gamma_{f,n}\ \text{is not a graph for some }1\le n\le M\bigr\},
\]
for which $K_{f,M}(\lb,\om)\searrow K_f(\lb,\om)$ as $M\to\infty$ (compare Definition~\ref{def:K_f}); it is decided by finitely many iterates of $\Gamma_{f,0}$ alone and is evaluated by the algorithm of Appendix~\ref{app:error}.

Proposition~\ref{prop:Kf-bound} also shows that $\sup_fK_f(\lb,\om)$ is a lower bound for $\eps_*(\lb,\om)$ and can be approached by choosing suitable initial graphs $f$; we therefore take the five graphs $f_1,\dots,f_5$ of Appendix~\ref{app:error} (see Fig.~\ref{fig:omega-epsilon-bound}) and set $K_{1200}^{\text{max}}(\om)=\max_iK_{f_i,1200}(\lb,\om)$, which agrees with $\max_iK_{f_i}(\lb,\om)$ to within $10^{-5}$; scanning $\eps$ upward from $K_{1200}^{\text{max}}(\om)$ until the computed invariant circle becomes non-graphic gives $\eps_*^{\text{up}}(\om)\ge\eps_*(\lb,\om)$. Thus
\[
\max_iK_{f_i}(\lb,\om)\ \le\ \eps_*(\lb,\om)\ \le\ \eps_*^{\text{up}}(\om).
\]
The two ends of the band have different accuracy: $K_{1200}^{\max}$ approximates the lower bound to $10^{-5}$, while $\eps_*^{\mathrm{up}}$, obtained by scanning in steps of $10^{-3}$, is the value at which the graph first becomes non-graphic. For the sampled $\om$ the difference between the two curves is less than $10^{-3}$ in most cases, and for a few $\om$ it reaches $10^{-2}$, at most $1.5\times10^{-2}$; hence the two curves of Fig.~\ref{fig:omega-epsilon-bound} coincide over most of $[0,1/2]$. This also indicates that the true curve $\eps_*(\om)$ must lie between them. Proposition~\ref{prop:Kf-bound} does not assert $\sup_fK_f=\eps_*$: the supremum of these computable lower bounds is not known to reach the threshold, and enlarging the family $f_1,\dots,f_5$ can only improve it.

One can see that the curve $\eps_*(\om)$ shows pronounced dips in some regions. This is because, for the DSM, there is a structure analogous to the Arnold tongues of \cite{Arn}: if the $p/q$ Arnold tongue is defined by the existence of $p/q$-periodic points, then inside the tongue there is a region in which the existing periodic orbits of period $p/q$ are foci. A focus necessarily forces the invariant graph to become non-graphic, so $\eps_*(\om)$ avoids, or stays just below, the regions where such foci exist; these dips arise in this way.

Although numerically $\eps_*(\lb,\om)$ looks continuous, or even fairly smooth, we cannot be sure whether it really is continuous. Reference~\cite{KM} shows an example with a fractal critical boundary; see also~\cite{Si}. Compared with the area-preserving case, the dynamics of the dissipative standard map is more complicated, which makes it harder to verify the continuity of $\eps_*(\lb,\om)$. Therefore we turn to analysing the mechanism causing the breakdown.

\subsection{Lock of drift at $\eps_*(\lb,\om)$}

\begin{figure}[h]
\begin{center}
\includegraphics[width=13cm]{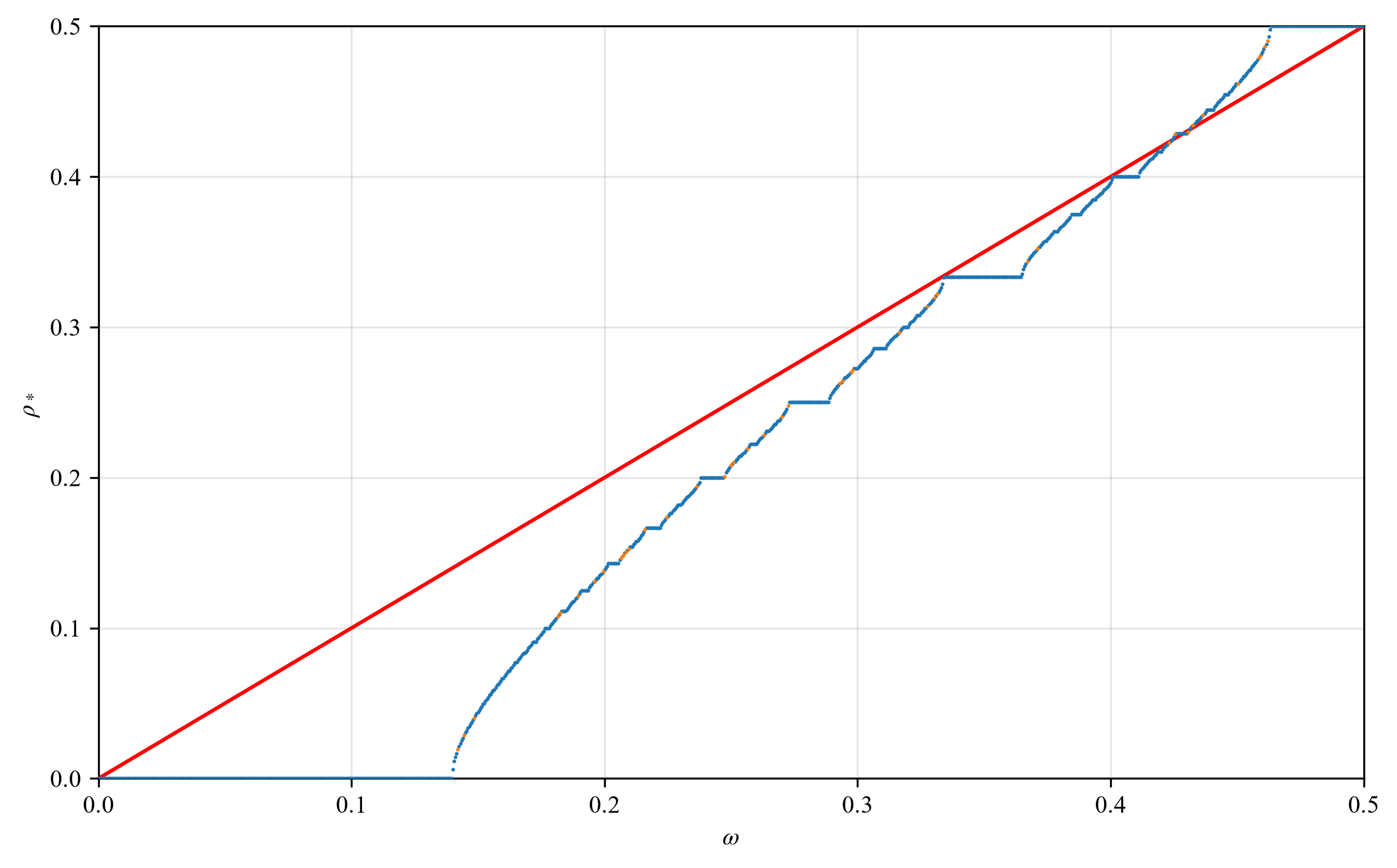}
\caption{Fix $\lb=2$. For 1000 sample values of $\omega$ equidistantly chosen in $[0,0.5]$, the discrete points give $\lim_{\eps\to\eps_*(\omega)}\rho(\omega,\eps)$. The values at the blue points are all exact rational numbers, whereas those at the orange points could not be determined. Compared with the diagonal (red line), the drift is not one-sided: its sign changes near $\om\approx0.42$. For $\om\le0.134$, $\lim_{\eps\to\eps_*(\omega)}\rho(\omega,\eps)$ is tightly locked at zero.}
\label{fig:omega-rho-501}
\end{center}
\end{figure}

\begin{figure}[h]
  \centering
  \begin{subfigure}[b]{0.48\textwidth}
    \centering
    \includegraphics[width=\linewidth]{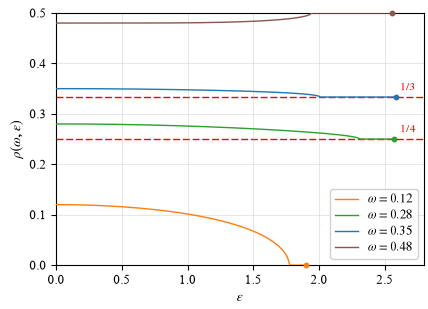}
    \caption{}
    \label{fig:epsilon-rho-4}
  \end{subfigure}
  \hfill
  \begin{subfigure}[b]{0.48\textwidth}
    \centering
    \includegraphics[width=\linewidth]{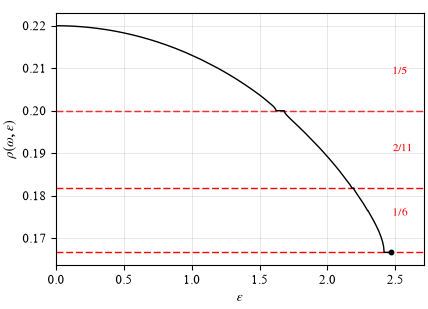}
    \caption{}
    \label{fig:anihilation}
  \end{subfigure}
  \caption{(A) Fix $\lambda=2$. This figure plots $\rho$ as a function of $\epsilon$ for four different values of $\omega$, along with relevant reference lines. The values of $\om$ and the other parameters are marked in the figure. Each curve terminates at the $\epsilon_*$ value corresponding to its $\omega$. (B) $\rho(\lambda,\om,\eps)$ for $\lambda=2$ and $\omega=0.22$.}
  \label{fig:epsilon-rho}
\end{figure}

\begin{figure}[h]
  \centering
  \begin{subfigure}[b]{0.48\textwidth}
    \centering
    \includegraphics[width=\linewidth]{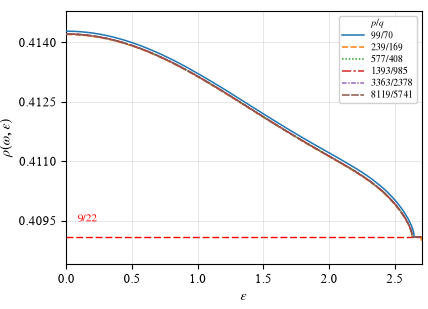}
  \end{subfigure}
  \hfill
  \begin{subfigure}[b]{0.48\textwidth}
    \centering
    \includegraphics[width=\linewidth]{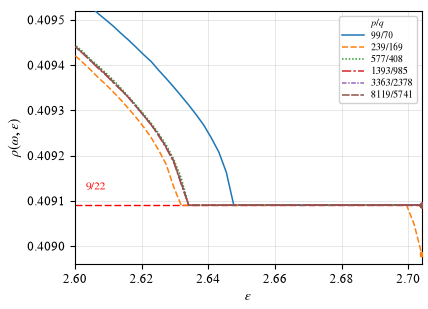}
  \end{subfigure}
  \caption{$\lambda=2$. We study $\om_*=\sqrt2-1\approx0.41421356$, approximated by six convergents of $\sqrt2$: $\om=p/q-1$ with $p/q=99/70$, $239/169$, $577/408$, $1393/985$, $3363/2378$ and $8119/5741$. (A) The rotation number $\rho(\om,\eps)$ for $\eps\in[0,2.704]$: all six curves start at their $\om$, decrease as $\eps$ grows, and lock onto the same rational value $9/22$ (red dashed line). (B) Zoom near the breakdown: the locking happens at slightly different values of $\eps$; the curve $p/q=239/169$ leaves $9/22$ at $\eps\approx2.70$, whereas the finer approximants $1393/985$, $3363/2378$ and $8119/5741$ keep the value $9/22$ up to the end of the window. The limit $9/22$ is therefore stable under refinement.}
  \label{fig:epsilon-rho-sqrt2}
\end{figure}

For the rotation-number drift, we fix $\lambda>0$ and a chosen $\omega$ and follow $\rho(\lambda,\omega,\epsilon)$ as $\epsilon$ increases toward its threshold. The numerical plots use $\lambda=2$: Figure~\ref{fig:omega-rho-501} compares the limiting rotation numbers across different $\omega$, whereas Figures~\ref{fig:epsilon-rho-4} and~\ref{fig:epsilon-rho-sqrt2} follow $\epsilon$-dependent curves for selected fixed values of $\omega$.

The variational observations at the beginning of this subsection hold for general $V\in C^2(\T,\R)$ (equivalently $V'\in C^1$). The rotation-locking evidence, figures, and the question below use the particular choice $V'(x)=\sin(2\pi x)$; the numerical claims are not asserted for every potential.

For such a potential, as $\eps>0$ increases while remaining small, the homogeneous invariant graph $\cT_{\eps=0}:=\{(x,0)\in T^*\T\}$ gets deformed slightly and maintains the graphic structure. By Theorem~\ref{thm:aubry}, $\cT_\eps$ coincides with $\wt\cA$. As an invariant graph, $\Phi:\wt\cA\rightarrow\wt\cA$ acts as an orientation-preserving circle homeomorphism, so its rotation number is independent of the orbit; we denote it by $\rho(\lb,\om,\eps)$. In the numerical experiments reported below, $\rho(\lb,\om,\eps)$ varies continuously with $\eps$ for as long as the invariant graph persists; together with $\rho(\lb,\om,0)=\om$, this allows the drift of the rotation number to be followed up to the threshold $\eps_*(\lb,\om)$.

Numerically, for many rational $\om=p/q$ there exists $\eps_-<\eps_*$ such that $\rho$ is constant on $(\eps_-,\eps_*)$ (Fig.~\ref{fig:epsilon-rho-4}). Such locking indicates that the invariant circle is then a stable structure formed by periodic orbits and their heteroclinic connections. Every nondegenerate periodic orbit is of one of three types: node, focus or saddle; as $\eps$ varies, the attracting periodic point may undergo a bifurcation from one type to another, while the heteroclinic manifolds also deform with $\eps$. Note that $\rho(\lb,\om,\eps_-)$ need not equal $\om$, even when $\om=p/q$ is rational with small denominator (Fig.~\ref{fig:omega-rho-501}).

It should be noted, however, that not all $p/q$-periodic points (i.e. solutions of $\Phi^q(x,y)=(x+p,y)$) lie on the invariant circle, but the points lying off the circle do not affect its rotation number (see Section~\ref{subsec:sga}). We focus on the case in which $\eps_-$ exists; in this case, the invariant circle, at the critical moment when it changes from a graph to a non-graphic curve, must contain $p/q$-periodic points. These periodic points therefore yield two breakdown mechanisms, treated in the next two subsections.

Another fact we want to point out is that, when $\om$ is irrational, the difference $(\eps_*-\eps_-)$ does not necessarily vanish. In some sense, even for irrational $\om\in[0,1/2]$, the lock of the drift still exists (the case $\om=\sqrt2-1$ is shown in Fig.~\ref{fig:epsilon-rho-sqrt2}). That makes it possible to clarify all the heteroclinic catastrophes leading to the breakdown.

\begin{que}
For $V'(x)=\sin(2\pi x)$ and fixed $\lb>0$, is there a uniform $\hat \eps_->0$ such that the associated $\rho(\lb,\om,\eps)$ remains a rational constant for all $\om\in[0,1/2]$ and $\eps\in[\eps_*-\hat\eps_-,\eps_*)$?
\end{que}

Locking may occur repeatedly without breakdown: Figure~\ref{fig:epsilon-rho}(B) shows $\rho$ locked at $1/5$ and $2/11$, then drifting while the graph persists, until the final locking at $1/6$. Thus rational locking alone does not imply breakdown.

\subsection{Saddle-node bridge}\label{subsec:saddle-node}

Throughout this subsection, $\lambda$ and $\om$ are fixed and $\eps$ varies. Motivated by the prevalent rational locking discussed in the previous subsection, we assume that the invariant circle consists of heteroclinic orbits connecting periodic points. The mechanism described in this subsection shows that the deformation of these heteroclinic orbits is one of the causes of the loss of the graphic property. The linear criterion below is independent of $V$; applying it to the DSM requires $V\in C^\infty(\T,\R)$ and a node with distinct positive eigenvalues, with resonances discussed in Lemma~\ref{lem:stowe} and Remark~\ref{rmk:resonance}. The numerics use $V'(x)=\sin(2\pi x)$, $\lambda=2$, $\om=0.1$.

\begin{figure}[htbp]
  \centering

  \begin{subfigure}[b]{0.48\textwidth}
    \centering
    \includegraphics[width=\linewidth]{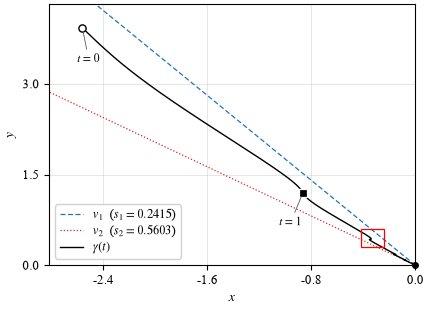}
    \caption{}
    \label{fig:model-1-a}
  \end{subfigure}
  \hfill
  \begin{subfigure}[b]{0.48\textwidth}
    \centering
    \includegraphics[width=\linewidth]{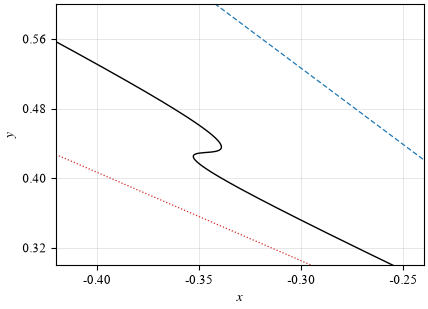}
    \caption{}
    \label{fig:model-1-b}
  \end{subfigure}

  \caption{The figure shows the linear map $z\mapsto Az$; $A=D\Phi(z_0)$ at the attracting node of the DSM (\ref{eq:stan}), $V'(x)=\sin 2\pi x$, $\lambda=2$, $\omega=0.1$, $\epsilon=1.734$, eigenvalue $0<s_1<s_2<1$, eigendirections $v_1,v_2$. Solid: $\gamma(t)=P\Lambda^{t-1}(1,1+0.25\sin 2\pi t)$, $t\in[0,\infty)$, $A=P\Lambda P^{-1}$, $\Lambda=\operatorname{diag}(s_1,s_2)$, $P=(v_1\ v_2)$. Node at the origin; (b) magnifies the red box, where $\gamma$ folds and is not a graph.}
  \label{fig:model-1}

\end{figure}

To understand the intrinsic mechanism, we first consider a linear model. Let $A$ be a linear map of $\R^2$ having the origin as an attracting node. The branch of the invariant manifold entering the origin is generated by forward iteration of a fundamental segment: any curve $\Sigma_0$ joining a point $P_0$ to $P_1:=AP_0$ gives the forward-invariant branch
\begin{equation}
  \Sigma^+ := \bigcup_{n\ge0}A^n(\Sigma_0),
\end{equation}
which accumulates at the origin. So we can make the following criterion for $\Sigma^+$:

\begin{prop}\label{prop:crit}
  Suppose $A$ has eigenvalues $0 < s_1 < s_2 < 1$ and neither eigendirection is vertical, and let $v_s$, $v_w$ be the eigenvectors of $s_1$, $s_2$ with positive $x$-components. Assume $\Sigma_0$ is an $x$-graph ($x$ is strictly monotone along $\Sigma_0$), with $0 \notin \Sigma_0$ and $x_0 := x(P_0) \neq 0$. Then $\Sigma^+$ is an $x$-graph if and only if, at every point of $\Sigma_0$, the oriented tangent vector $\gamma'(t) = \alpha(t)\, v_s + \beta(t)\, v_w$ satisfies
  \begin{equation*}
    \beta(t)\, x(\gamma'(t)) \ge 0.
  \end{equation*}
  In other words, the admissible sector for the oriented tangent is determined by $A$ alone; once the oriented tangent leaves this sector, the forward iterates necessarily become non-graphic.
\end{prop}

\begin{proof}
  Parametrize $\Sigma_0$ by $\gamma \colon [0,1] \to \mathbb{R}^2$, $\gamma(0) = P_0$, $\gamma(1) = P_1$; since $\Sigma_0$ is an $x$-graph, $c(t) := x(\gamma'(t))$ has constant sign, vanishing only at isolated vertical tangencies. Reflecting in the $y$-axis if necessary (which preserves the condition $\beta \cdot x(\gamma') \ge 0$ and interchanges the two monotonicity directions), assume $x_0 > 0$.

  For any direction $d = \alpha v_s + \beta v_w$,
  \begin{equation*}
    A^n d = s_2^n \big( \beta v_w + \alpha t_n v_s \big), \qquad
    t_n := (s_1/s_2)^n \in (0,1].
  \end{equation*}
  If $\beta \ge 0$ and $x(d) \ge 0$, then $x(A^n d) = s_2^n(\beta x_w + \alpha x_s t_n) \ge 0$, strictly for $n \ge 1$ and $d \neq 0$, because the linear function $\beta x_w + \alpha x_s t$ has endpoint values $\beta x_w \ge 0$ and $x(d) \ge 0$ and vanishes identically only for $d = 0$. Hence the quadrant $\{\beta \ge 0,\ x \ge 0\}$ is forward invariant; likewise $\{\beta \le 0,\ x \le 0\}$ by the same argument.

  \emph{Sufficiency.} Suppose $\beta \cdot c \ge 0$ everywhere. If $c \ge 0$, then $\gamma'$ lies in $\{\beta \ge 0,\ x \ge 0\}$, so with $h_n(t) := x(A^n \gamma'(t))$ we get $x_{n+1} - x_n = \int_0^1 h_n dt \ge 0$ and $x_n \ge x_0 > 0$ for all $n$, contradicting $x_n = x(A^n P_0) \to 0$. Hence $c \le 0$ and $\gamma'$ lies in $\{\beta \le 0,\ x \le 0\}$; then $h_n < 0$ a.e.\ for $n \ge 1$, so $x \circ A^n \gamma$ is strictly decreasing and $x_{n+1} < x_n$, $x_n \downarrow 0$. Each $A^n \Sigma_0$ is therefore an $x$-graph over $[x_{n+1}, x_n]$; these intervals tile $(0, x_0]$, so $\Sigma^+$ is an $x$-graph.

  \emph{Necessity.} Suppose $\Sigma^+$ is an $x$-graph, say $\Sigma^+ = \{(x, f(x)) : 0 < x \le x_0\}$. The map $x' := \pi_x \circ A \circ (\mathrm{id}, f)$ is continuous and fixed-point free on $(0, x_0]$: a fixed point would be a fixed point of $A$ lying on $\Sigma^+$, while the only fixed point of $A$ is $0 \notin \Sigma^+$. Hence $x' - \mathrm{id}$ has constant sign, and $x_n = (x')^n(x_0) \to 0$ forces $x' < \mathrm{id}$. Thus $x_n \downarrow 0$, each $A^n \Sigma_0$ is the sub-arc of the graph over $[x_{n+1}, x_n]$, and along $\gamma$ the coordinate $x$ decreases, so $h_n \le 0$ and in particular $c \le 0$. If $\beta(t_0) > 0$ at some point $t_0$, then $\beta > 0$, $c < 0$ and hence $\alpha < 0$ (from $c = \alpha x_s + \beta x_w < 0$) on a neighbourhood $U$ of $t_0$. Since $\beta x_w + \alpha x_s t_n \to \beta x_w > 0$ uniformly on $U$, we get
  \begin{equation*}
    h_n(t) = s_2^n \big( \beta x_w + \alpha x_s t_n \big) > 0 \qquad
    \text{on } U
  \end{equation*}
  for $n$ large, contradicting $h_n \le 0$. Hence $\beta \le 0$ everywhere, and $\beta \cdot c \ge 0$.
\end{proof}

\begin{rmk}
We have assumed that neither eigendirection is vertical. This entails no loss of generality for the systems considered here: for the node or saddle of the DSM $\Phi_{\lambda,\om,\epsilon}$ lying on the invariant circle, its eigendirections cannot be vertical, as proved in Corollary~\ref{cor:graph-periodic-linear} of Section~\ref{s2}; hence the vertical case need not be discussed separately, and the two cone boundaries in Proposition~\ref{prop:crit} are well-defined. Consequently, in the saddle--node bridge scenario the loss of the graphic property cannot be attributed to an eigendirection becoming vertical while the orbit remains globally minimizing; this does not exclude other breakdown mechanisms, such as the node-to-focus transition of Subsection~\ref{subsec:bifurcation}.
\end{rmk}

For the DSM, let $z_0$ be an attracting $p/q$-periodic node, i.e.\ a point with $\Phi^{q}(z_0)=z_0+(p,0)$ whose differential
\begin{equation}\label{eq:Aq}
  A:=D\Phi^{q}(z_0)
\end{equation}
has distinct eigenvalues $0<s_1<s_2<1$; for $q=1$ this is a fixed node, and the figures below concern this case. By Corollary~\ref{cor:graph-periodic-linear} neither eigendirection of $A$ is vertical. Put
\begin{equation}\label{eq:spread}
  \sigma:=\frac{\log s_1}{\log s_2}\in(1,+\infty);
\end{equation}
note that $\sigma\to1$ as the eigenvalues collide at the node--focus transition. Since $s_1s_2=\det A=e^{-\lb q}$ (Corollary~\ref{cor:graph-periodic-linear}), the equality $s_1=s_2^{m}$ with an integer $m\ge2$ is equivalent to $\sigma=m$; these are the only possible resonances of $A$, at any order.

The following sharp statement links the order of a resonance to the smoothness of a linearizing change of variables.

\begin{lem}[Linearization, {\rm\cite[Theorems~4--6]{St}}]\label{lem:stowe}
Assume $V\in C^\infty(\T,\R)$, so that $\Phi^{q}\in C^\infty$, and let $z_0$ be a $p/q$-periodic point for which $A=D\Phi^{q}(z_0)$ has distinct real eigenvalues $0<s_1<s_2<1$.
\begin{itemize}
\item[(i)] If $\sigma\notin\Z$, there are a neighbourhood $U$ of $z_0$ and a $C^\infty$ diffeomorphism $\Psi:U\to\Psi(U)$ with $\Psi(z_0)=0$, $D\Psi(z_0)=\mathrm{id}$ and $\Psi\circ\Phi^{q}\circ\Psi^{-1}=A$ on $\Psi(U)$.
\item[(ii)] If $\sigma=m\in\Z$, the same conclusion holds with $\Psi$ of class $\bigcap_{\varepsilon>0}C^{m-\varepsilon}$, and in general not of class $C^{m}$.
\end{itemize}
\end{lem}

Thus for a $C^\infty$ map the resonance order is exactly the number of derivatives that the linearization loses: none if $\sigma\notin\Z$, and $m$ if $\sigma=m$. Since the possible orders are $m=2,3,\dots$, only the lowest one can destroy $C^2$ regularity: for $m\ge3$ the class $C^{m-\varepsilon}$ still contains $C^2$, and only $m=2$ gives $C^{2-\varepsilon}$. (If one assumes merely $V\in C^3$, so that $\Phi^{q}\in C^2$, the same circle of ideas yields $C^2$ linearization exactly for $\sigma<2$; the $C^\infty$ hypothesis is what allows the resonance orders $m\ge3$ to be absorbed.)

The following corollary shows that the cone criterion from the linear model can be transferred to the nonlinear system with an explicitly estimable error.

\begin{cor}\label{cor:cone-error}
Assume $V\in C^\infty(\T,\R)$ and $\sigma\ne2$, and let $\Psi$ be a $C^2$ linearization as in Lemma~\ref{lem:stowe}. For $\eta>0$ put
\begin{equation}\label{eq:Meta}
  M_\eta:=\sup_{B(z_0,\eta)}\|D^2\Psi\|<+\infty,\qquad \kappa:=M_\eta\eta .
\end{equation}
Then for $z\in B(z_0,\eta)$ and every direction $d$ the transported direction $\widehat d:=D\Psi(z)d$ satisfies
\begin{equation}\label{eq:angle-error}
  \angle(d,\widehat d)\le\arcsin\frac{\kappa}{1-\kappa}\qquad(\kappa<1).
\end{equation}
\end{cor}
\begin{proof}
By the mean value theorem $\|D\Psi(z)-\mathrm{id}\|\le M_\eta|z|\le\kappa$, uniformly for $z\in B(z_0,\eta)$. Hence the component of $\widehat d$ perpendicular to $d$ is at most $\kappa|d|$, while $|\widehat d|\ge(1-\kappa)|d|$.
\end{proof}

\begin{figure}[h]
  \centering

  \begin{subfigure}[b]{0.48\textwidth}
    \centering
    \includegraphics[width=\linewidth]{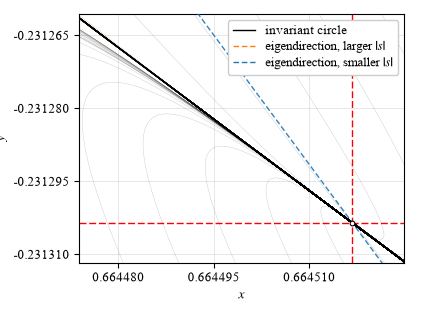}
    \caption{}
    \label{fig:1pp-example-2}
  \end{subfigure}
  \hfill
  \begin{subfigure}[b]{0.48\textwidth}
    \centering
    \includegraphics[width=\linewidth]{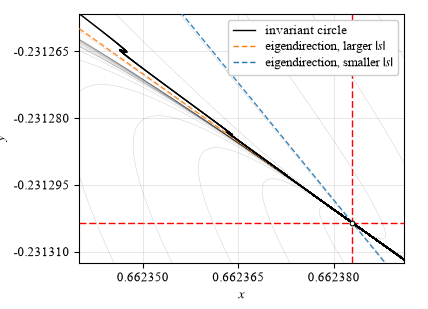}
    \caption{}
    \label{fig:1pp-example-4}
  \end{subfigure}

  \caption{$\lambda=2$, $\omega=0.1$. Left: $\eps=1.72$, the invariant circle is still graphic in the zoomed window. Right: $\eps=1.734$; it has become non-graphic. The horizontal and vertical axes represent the coordinates $(x,y)$ on the cylinder, and the windows are indicated in the subfigures. The thin gray lines are the flow of the corresponding continuous dynamical system.}
  \label{fig:1pp-example}
\end{figure}

\begin{rmk}\label{rmk:resonance}
Two points deserve emphasis.

\emph{(1) The resonant set is codimension one.} By \eqref{eq:spread} and $s_1s_2=e^{-\lb q}$, the condition $\sigma=m$ is the single scalar equation $s_1=s_2^{m}$, equivalently $s_2=e^{-\lb q/(m+1)}$; the lowest member, the second-order resonance, reads
\begin{equation}\label{eq:sigma2}
  s_1=s_2^2\qquad\Longleftrightarrow\qquad s_2=e^{-\lb q/3}.
\end{equation}
Since $\sigma$ is a continuous function of the parameters, the set where it takes integer values is a countable union of level sets of one function; wherever $\nabla\sigma\ne0$ these are codimension-one submanifolds. Hence resonances are expected only on a codimension-one part of parameter space, and along a one-parameter family they occur at isolated parameters, transversally whenever $\partial_\eps s_2\ne0$. In particular $C^\infty$ linearization, and a fortiori $C^2$ linearization, holds for all parameters outside this set, i.e.\ generically; for instance the nodes of Fig.~\ref{fig:1pp-example} have $\sigma=2.64$ and $2.45$, neither of them an integer.

\emph{(2) Only the second-order resonance can spoil $C^2$.} This is Lemma~\ref{lem:stowe} combined with the dictionary above: at $\sigma=m\ge3$ the linearization is still $C^{m-\varepsilon}\subset C^2$, whereas at $\sigma=2$ it is only $C^{2-\varepsilon}$, so that the constant $M_\eta$ in \eqref{eq:Meta} is unavailable and the cone can be compared with the linear model only up to an angular error $O(\eta^{1-\varepsilon})$. Note also that $\sigma\to1$ as the node approaches the node--focus transition, so along such a branch the value $\sigma=2$ is met at most at isolated parameters.
\end{rmk}

\begin{figure}[htbp]
  \centering

  \begin{subfigure}[b]{0.48\textwidth}
    \centering
    \includegraphics[height=5.0cm]{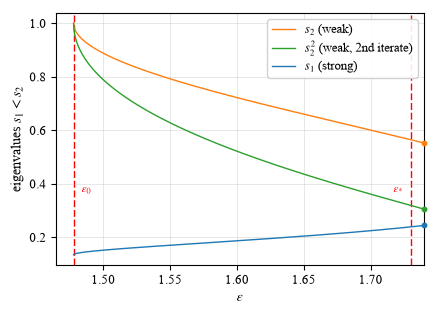}
    \caption{}
    \label{fig:eigenvalue-1pp}
  \end{subfigure}
  \hfill
  \begin{subfigure}[b]{0.48\textwidth}
    \centering
    \includegraphics[height=5.0cm]{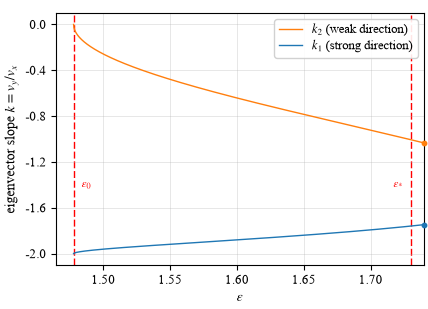}
    \caption{}
    \label{fig:eigendirect-1pp}
  \end{subfigure}

  \caption{$\lambda=2$, $\omega=0.1$. (A) The multipliers $s_1<s_2$ along the node branch, together with the powers $s_2^2$: the ratio $\sigma=\log s_1/\log s_2$ equals $2$ at the resonance $s_1=s_2^2$, but this does not occur before the invariant circle becomes non-graphic. (B) Both eigendirection slopes stay finite. Dashed lines mark the appearance of the node and $\epsilon_*$.}
  \label{fig:eigen_val_dir}

\end{figure}

By Lemma~\ref{lem:stowe}, the admissible sector of the linear model is carried to a curved sector $\Psi(\cD)$ near the node, where
\begin{equation*}
  \cD := \left\{ \, d = \alpha v_s + \beta v_w \;:\; \beta \cdot x(d) \ge 0 \, \right\}.
\end{equation*}
The sector is invariant, $A(\cD) \subset \cD$ (cf.\ the proof of Proposition~\ref{prop:crit}); its strong boundary $\operatorname{span}\{v_s\}$ is carried by $\Psi$ to an invariant curve tangent to the strong eigendirection at the node. Since $D\Psi(z_0) = \mathrm{id}$, the curved sector $\Psi(\cD)$ and the linear sector $\cD$ have the same tangent cone at the node; quantitatively, on $B(z_0,\eta)$ their angular separation is at most $\arcsin(\kappa/(1-\kappa))$ with $\kappa=M_\eta\eta$ whenever $\sigma\ne2$ (Corollary~\ref{cor:cone-error}), and at most $O(\eta^{1-\varepsilon})$ at $\sigma=2$ (Remark~\ref{rmk:resonance}). Proposition~\ref{prop:crit} therefore applies near the node with $\cD$ replaced by $\Psi(\cD)$, up to this angular error.

As $\epsilon \to \epsilon^*$, the invariant circle, being the union of unstable branches of the saddle-type periodic points each accumulating on a node, is stretched vertically, so each branch enters the node with an increasingly extreme slope. Then its oriented tangent vector may cross the admissible sector determined by $D\Phi_{\lb,\om,\eps}$; by Proposition~\ref{prop:crit}, the forward iterates of the unstable manifold then become non-graphic, exactly as in the linear model, up to the angular error of Corollary~\ref{cor:cone-error}. This is the mechanism shown in Fig.~\ref{fig:1pp-example}.

\subsection{Bifurcation of node points}\label{subsec:bifurcation}

\begin{figure}[h]
  \centering

  \begin{subfigure}[b]{0.32\textwidth}
    \centering
    \includegraphics[width=\linewidth]{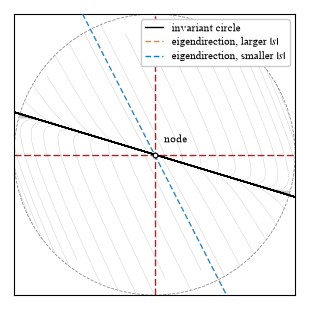}
    \caption{}
    \label{fig:Omega0-AL}
  \end{subfigure}
  \hfill
  \begin{subfigure}[b]{0.32\textwidth}
    \centering
    \includegraphics[width=\linewidth]{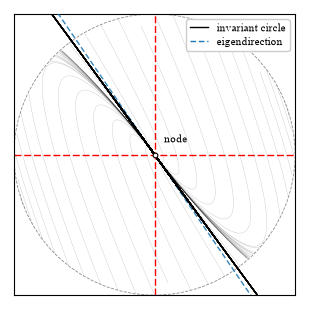}
    \caption{}
    \label{fig:Omega0-BL}
  \end{subfigure}
  \hfill
  \begin{subfigure}[b]{0.32\textwidth}
    \centering
    \includegraphics[width=\linewidth]{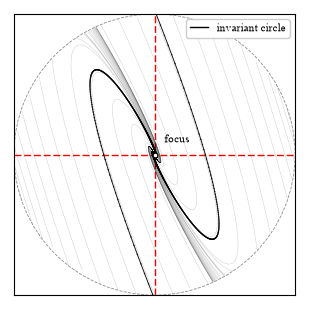}
    \caption{}
    \label{fig:Omega0-CL}
  \end{subfigure}

  \caption{With $\lambda=2$ fixed, the three panels show the local picture near the attractor for $\epsilon=0.3$, $\eps_*(2,0)$ and $1.9$, respectively; the corresponding attractors are a node, an improper node and a focus. One can see that in the first two cases the invariant circle is graphic. The radius of the circles in the panels is $0.0005$.}
  \label{fig:Omega0}
\end{figure}

The second mechanism is the bifurcation of  nodes. As $\eps$ rises, the two real eigenvalues of an attracting node may collide and then become a complex conjugate pair, so that the node turns into a focus; since a focus attracts spirally, the invariant graph cannot survive such a bifurcation. Here $\lambda>0$ and $\om$ are fixed while $\eps$ varies. The explicit calculation takes $\om=0$ (Fig.~\ref{fig:Omega0} specializes to $\lambda=2$) and $V'(x)=\sin(2\pi x)$.

As an example, for $\om=0$ this mechanism can be analysed completely. Fix $\lb$; for every $\eps$ the points $(0,0)$ and $(1/2,0)$ are fixed points of \eqref{eq:stan}. The former is a saddle for all $\eps$, while the latter is an attracting node for small $\eps$: as $\eps$ increases, its two real eigenvalues approach each other and coincide at $\eps=\eps_c$, where the node is improper, and for $\eps>\eps_c$ they form a complex conjugate pair, so that the point becomes a focus. A direct computation gives
\[
\eps_c=\frac{(e^{\lb/2}-1)^2}{2\pi\mu}=\frac{\lb e^{\lb}}{2\pi}\tanh\frac{\lb}{4}=K_0(\lb,0),
\]
the last identity being Corollary~\ref{cor:K0zero} for $V'(x)=\sin2\pi x$.

By Theorem~\ref{thm:e_u_i} applied to the zero section, for every $\eps<K_0(\lb,0)=\eps_c$ there is a unique invariant graph, so $\eps_*(\lb,0)\ge\eps_c$. But for $\eps>\eps_c$ the fixed point $(1/2,0)$ is a focus, and no invariant circle passes through a focus. If an invariant graph existed, it would be the unique invariant circle, since the uniqueness argument in the proof of Theorem~\ref{thm:e_u_i} does not use $\eps<K_f$. The map with $\om=0$ commutes with $S(x,y)=(1-x,-y)$, so this circle would be $S$-invariant; a graph $\{(x,\psi(x))\}$ invariant under $S$ satisfies $\psi(1-x)=-\psi(x)$, hence $\psi(1/2)=0$ and the graph passes through $(1/2,0)$, a contradiction. Therefore no invariant graph exists for $\eps>\eps_c$, and $\eps_*(\lb,0)\le\eps_c$. Consequently
\[
\eps_*(\lb,0)=\eps_c=K_0(\lb,0)=\frac{\lb e^{\lb}}{2\pi}\tanh\frac{\lb}{4}.
\]
At $\om=0$ the zero section therefore already reaches the threshold: $K_0(\lb,0)=\eps_*(\lb,0)$.

See Fig.~\ref{fig:Omega0}.

When $\omega=1/2$, a similar mechanism can also be obtained by direct computation. Beyond this, there are other values of $\omega$ for which this mechanism is responsible for the breakdown of the invariant graph, although these can no longer be obtained in closed form.

\begin{rmk}\label{rmk:KMW-breakdown}
For $\omega=0$, \cite{KMW} also shows that $\eps>K_0(\lambda,0)$ is a sufficient condition for the breakdown of the invariant graph.
\end{rmk}

\subsection{Invariant-graph breakdown and the global attractor}\label{subsec:sga}

The breakdown of an invariant graph should be distinguished from the development of a more complicated global attractor. In particular, the rotation number of the invariant graph need not coincide with either of the two rotation numbers associated with the boundary of the global attractor.

Following Le Calvez~\cite{Cal}, the upper and lower boundaries of the global attractor carry rotation numbers $\rho^+$ and $\rho^-$, respectively.

\begin{defn}[Surfaces and boundary rotation numbers, Le Calvez \cite{Cal}]\label{def:surfaces}
Let $\widetilde{\Phi}:\mathbb{R}^{2}\to\mathbb{R}^{2}$ be the lift of $\Phi=\Phi_{\lambda,\omega,\epsilon}$ satisfying $\widetilde{\Phi}(x+1,y)=\widetilde{\Phi}(x,y)+(1,0)$, and let $\pi_{1}:\mathbb{R}^{2}\to\mathbb{R}$ be the angular projection. Recall the global attractor $\Omega$ of Definition~\ref{def:ga}. For each $x\in\mathbb{T}$ set
\[
y^{-}(x):=\min\{y\mid(x,y)\in\Omega\},\qquad
y^{+}(x):=\max\{y\mid(x,y)\in\Omega\};
\]
their graphs $\partial\Omega^{\pm}:=\{(x,y^{\pm}(x))\mid x\in\T\}$ are called the {\bf lower/upper surfaces} of $\Omega$ ($y^{-}$ is lower, and $y^{+}$ upper, semi-continuous). Since $\widetilde{\Phi}^{-1}(\partial\Omega^{\pm})\subset\partial\Omega^{\pm}$, the maps
\[
G^{\pm}(x):=\pi_{1}\bigl(\widetilde{\Phi}^{-1}(x,y^{\pm}(x))\bigr)
\]
are non-decreasing and satisfy $G^{\pm}(x+1)=G^{\pm}(x)+1$, hence the {\bf lower/upper rotation numbers}
\[
\rho^{\pm}(\lambda,\omega,\epsilon)
:=-\lim_{n\to+\infty}\frac{(G^{\pm})^{n}(x)-x}{n}
\]
are well-defined and independent of $x$ \cite[\S 4--5]{Cal}.
\end{defn}

The Aubry set $\wt\cA$ is always contained in $\Omega$, and its rotation number lies between the two boundary rotation numbers. In particular, as long as the invariant graph $\cT_\epsilon$ exists, Theorem~\ref{thm:aubry} identifies it with $\wt\cA$, so that
\begin{equation}\label{eq:rho-surfaces}
\rho^{-}(\lambda,\omega,\epsilon)
\leq \rho(\lambda,\omega,\epsilon)
\leq \rho^{+}(\lambda,\omega,\epsilon).
\end{equation}
The inequalities in \eqref{eq:rho-surfaces} need not be equalities: the rotation number of the invariant graph can lie strictly between $\rho^-$ and $\rho^+$, or coincide with one or both of the boundary rotation numbers. Thus the dynamics on the boundary of the global attractor does not, in general, determine the rotation number of the invariant graph.

\begin{figure}[htbp]
  \centering

  \begin{subfigure}[b]{0.48\textwidth}
    \centering
    \includegraphics[width=\linewidth]{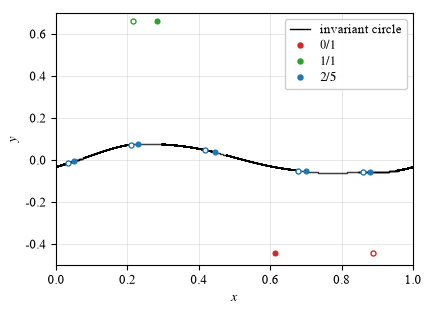}
    \caption{}
    \label{fig:bothside}
  \end{subfigure}
  \hfill
  \begin{subfigure}[b]{0.48\textwidth}
    \centering
    \includegraphics[width=\linewidth]{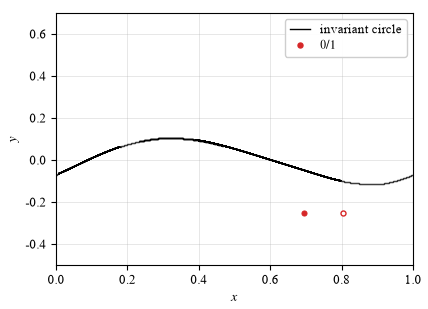}
    \caption{}
    \label{fig:oneside}
  \end{subfigure}

  \caption{Fix $V'(x)=\sin(2\pi x)$. Periodic points in the global attractor are marked, hollow symbols being saddles and solid symbols nodes. (A) $\lambda=0.2$, $\omega=0.4$, $\epsilon=0.15$: the invariant graph (red) has rotation number $2/5$, while the boundary rotation numbers of the attractor are $\rho^-=0$ and $\rho^+=1$, carried by the periodic orbits of rotation numbers $0/1$ and $1/1$. (B) $\lambda=0.5$, $\omega=0.2$, $\epsilon=0.175$: $\rho^-=0$ and $\rho^+=\rho$, so that the rotation number of the invariant graph coincides with the upper boundary rotation number.}
  \label{fig:attractor}
\end{figure}

Figure~\ref{fig:attractor} gives two illustrations. In panel (A),
\[
\rho^-=0<\rho=\frac{2}{5}<\rho^+=1,
\]
so the rotation number of the invariant graph differs from both boundary rotation numbers, which are carried by the periodic orbits of rotation numbers $0/1$ and $1/1$ visible in the attractor; in panel (B),
\[
\rho^-=0<\rho^+=\rho(\lb,\om,\eps),
\]
so that it coincides with one of them. In both panels the invariant graph is still graphic.

The breakdown of the invariant graph is not related to these periodic points. As $\eps$ increases, the rotation number of the invariant graph stays different from that of the $0/1$ orbit and from that of the $1/1$ orbit, and the loss of the graphic property occurs without any collision or resonance with them. Moreover, even after the graphic property is lost, the invariant circle continues to exist as an attractor and remains disjoint from these periodic orbits, which persist separately in the attractor: it is joined neither to the fixed point of rotation number $0$ nor to the $1/1$ orbit. Consequently the rotation number of the invariant graph is not constrained by $\rho^\pm$: it may coincide with both of them, with one of them, or with neither --- in panel (A) it lies strictly between them. Thus neither $\rho=\rho^\pm$ nor $\rho^-=\rho^+$ provides a general criterion for the breakdown threshold $\eps_*(\lb,\om)$; the mechanisms studied above concern the geometry of the invariant circle itself and should not be identified with changes in the global structure of $\Omega$.\medskip

\appendix
\section{Finite-iteration thresholds: truncation and computational error}\label{app:error}

\noindent{\bf The two thresholds.} For $f\in C(\T,\R)$ let $\Gamma_{f,0}=\{(x,f(x)):x\in\T\}$ and $\Gamma_{f,n}=\Phi^n(\Gamma_{f,0})$, and set
\[
K_{f,M}(\lb,\om):=\inf\bigl\{\eps\ge0:\ \Gamma_{f,n}\ \text{is not a graph for some }1\le n\le M\bigr\},
\]
with $K_f$ given by the same formula with $M$ replaced by $+\infty$ (Definition~\ref{def:K_f}). Then $K_f\le K_{f,M}$, the sequence $K_{f,M}$ is non-increasing in $M$, and $K_{f,M}\downarrow K_f$. The thresholds quoted in Section~\ref{s4} are $K_{f_i,1200}$, $i=1,\dots,5$, and their maximum, where
\[
f_1=0.20\sin2\pi x+0.1,\quad f_2=0.30\sin2\pi x+0.1,\quad f_3=0.30\sin2\pi x+0.2,
\]
\[
f_4=0.30\sin2\pi x+0.10\sin4\pi x,\quad f_5=0.20\sin2\pi x+0.10\sin4\pi x.
\]

\noindent{\bf Computation.} $\Gamma_{f,n}$ is a graph if and only if $x\mapsto X_n(x)$ is non-decreasing, i.e.\ $a_n(x)=\partial_xX_n(x)\ge0$ for all $x$, where $a_{n+1}=a_n+\mu b_n$ and $b_{n+1}=e^{-\lb}\eps V''(X_{n+1})a_n+e^{-\lb}(1+\eps\mu V''(X_{n+1}))b_n$, starting from $(a_0,b_0)=(1,f'(x))$. With $\hat a_n=a_n/\sqrt{a_n^2+b_n^2}$, the threshold $K_{f,M}$ is the first zero of $G_{f,M}(\eps)=\min_{1\le n\le M}\min_{x\in\T}\hat a_n(x;\eps)$, found by bisection in $\eps$. The minimum in $x$ is evaluated on $N$ equally spaced points and then refined by a local search; the resulting threshold is written $K^N_{f,M}$. Since a finite sample can only overestimate the minimum,
\begin{equation}\label{eq:onesided}
K_f\le K_{f,M}\le K^N_{f,M}.
\end{equation}

\noindent{\bf Truncation error.} Table~\ref{tab:trunc} lists $\max_iK_{f_i,M}$ for several $\om$ and $M$. For $\om$ away from $0$ and $1/2$ the values do not change from $M=300$ on, while at $\om=0$ and $\om=1/2$ they decrease like $M^{-2}$: at $M=1200$ the difference to $M=4800$ is $6.4\times10^{-6}$ at $\om=0$, $5.4\times10^{-6}$ at $\om=1/2$, and below $10^{-6}$ at the other sampled values. The values quoted in Section~\ref{s4} are therefore within $10^{-5}$ of $\max_iK_{f_i}$.

\begin{table}[htbp]
\centering
\caption{The truncation error: $\max_iK_{f_i,M}$ at $\lb=2$ for several $\om$ and $M$, computed with $N=2048$.}
\label{tab:trunc}
\begin{tabular}{cccccc}\hline
$\om$ & $M=300$ & $M=600$ & $M=1200$ & $M=2400$ & $M=4800$\\ \hline
0.00 & 1.0870133 & 1.0869314 & 1.0869108 & 1.0869057 & 1.0869044\\
0.10 & 1.7304725 & 1.7304725 & 1.7304725 & 1.7304725 & 1.7304725\\
0.25 & 2.5424152 & 2.5424152 & 2.5424152 & 2.5424152 & 2.5424152\\
0.40 & 2.6955775 & 2.6955775 & 2.6955775 & 2.6955775 & 2.6955775\\
0.50 & 2.3521025 & 2.3520329 & 2.3520154 & 2.3520111 & 2.3520100\\
\hline
\end{tabular}
\end{table}

\noindent{\bf Discretisation error.} Table~\ref{tab:grid} lists $\max_iK^N_{f_i,1200}$ for several grid sizes. The values decrease with $N$ and stabilise: at $\om=0.3$ and $\om=0.45$, where the dependence on $N$ is largest, the entries for $N=1024$ and $N=2048$ differ by $1.2\times10^{-5}$ and $2.0\times10^{-5}$, while the entries for $N=2048$ and $N=4096$ agree to $5\times10^{-8}$. The values quoted in Section~\ref{s4} were computed with $N=2048$; by \eqref{eq:onesided} they are upper bounds for $K_f$, with a discretisation error below $10^{-7}$.

\begin{table}[htbp]
\centering
\caption{The discretisation error: $\max_iK^N_{f_i,1200}$ at $\lb=2$ for several $\om$ and $N$, computed with $M=1200$.}
\label{tab:grid}
\begin{tabular}{cccc}\hline
$\om$ & $N=512$ & $N=1024$ & $N=2048$\\ \hline
0.00 & 1.0869108 & 1.0869108 & 1.0869108\\
0.10 & 1.7304828 & 1.7304814 & 1.7304725\\
0.30 & 2.6242502 & 2.6242502 & 2.6242384\\
0.45 & 2.6914873 & 2.6914604 & 2.6914409\\
\hline
\end{tabular}
\end{table}

\noindent Summing up, Tables~\ref{tab:trunc} and~\ref{tab:grid} show that the thresholds quoted in Section~\ref{s4} are accurate to $10^{-5}$.

\end{document}